\documentclass[a4paper,10pt]{amsart}
\usepackage{amssymb}

\numberwithin{equation}{section}
\newtheorem{theorem}{Theorem}[section]

\newtheorem{lem}[theorem]{Lemma}

\theoremstyle{remark}
\newtheorem{rem}[theorem]{Remark}
\newcommand{\R}{\mathbb{R}}
\newcommand{\Z}{\mathbb{Z}}
\newcommand{\N}{\mathbb{N}}

\renewcommand{\hat}{\widehat}

\author[J.~Benameur]{Jamel Benameur}
\author[S.~Ben Abdallah]{Saber Ben Abdallah}
\address{Department of Mathematics, Faculty of Science of Gab\`es, Research Laboratory Mathematics and Applications LR17ES11; Tunisia}
\email{\sl jamelbenameur@gmail.com}
\email{\sl saber\underline{\;\;}tradition@yahoo.fr}
\title[Asymptotic behavior study of critical ...]
{Asymptotic behavior of critical dissipative quasi-geostrophic equation in Fourier space}

\begin{document}
\begin{abstract}
In this paper we show the global existence for critical dissipative quasi-geostrophic equations if $\|\widehat{\theta^0}\|_{L^1}$ is small enough; among others we prove the analyticity of such a solution. If in addition the initial condition verifies $|D|^{-\delta}\theta^0\in L^1(\R^2)$ with $0<\delta<1$, then the solution remains regular and $\lim_{t\rightarrow\infty}t^\delta\|\widehat{\theta}(t)\|_{L^1}=0$. Fourier analysis and standard techniques are used.\end{abstract}


\subjclass[2000]{35-xx, 35Bxx, 35Jxx}
\keywords{Quasi-geostrophic Equations; Critical spaces; Long time decay}

\maketitle
\tableofcontents


\section{Introduction}
We interest to the study of dissipative quasi-geostrophic equations in some critical space represented by the following system:
$$ (S_1)\begin{cases}
  \partial_t\theta+\kappa|D|\theta+u_\theta.\nabla\theta  = 0\;\;{\rm in}\;\mathbb{R}^{+*}\times\mathbb{R}^2\\
u_\theta=\mathcal{R}^\bot\theta=(-R_2\theta,R_1\theta)\\
\theta_{/_{t=0}}=\theta^0,
   \end{cases}$$
where $R_1,R_2$ are the Riesz's transformation giving by $$\widehat{R_jf}=-\frac{i\xi_j}{|\xi|}\hat{f},$$
and $\kappa>0$ is the dissipative coefficient. for more details and mathematical and physical explanations of this model we can consult \cite{CLM}, \cite{CMT}, \cite{DC2} and \cite{JP}. These equations have the property of invariance by change of scales called "Scaling" as following: If $\theta(t,x)$  is solution of $(S_1)$ on $[0,T]$ with initial data $\theta^0=\theta(0,x)$ then $\theta_{\lambda}$ defined by $\theta_{\lambda}(t,x)=\theta(\lambda t,\lambda x)$ is a solution of $(S_1)$ on $[0,T/\lambda]$ with initial data $\theta_{\lambda}^0=\theta_{\lambda}(0,\lambda x)$. For $\sigma\in\R$, we define the de Fourier space by $${X}^{\sigma}(\mathbb{R}^2)=\{f\in\mathcal{S}'(\R^2);\;\hat{f}\in L_{loc}^1(\mathbb{R}^2)\;{\rm and}\;\int_{\mathbb{R}^2}|\xi|^\sigma|\hat{f}(\xi)|d\xi<\infty\},$$
which is equipped with the following norm $$\|f\|_{X^\sigma}=\int_{\mathbb{R}^2}|\xi|^\sigma|\hat{f}(\xi)|d\xi.$$
The critical space for $(S_1)$ is the Fourier's space definite previous for $\sigma=0$. We recall that a functional space $(X,\|.\|_X)$ is critical for the equations $(S_1)$ if $\|f(\lambda.)\|_X=\|f\|_X$ for all $\lambda>0$. Particularly, $L^\infty(\R^2)$, $\dot H^1(\R^2)$ and ${X}^0(\R^2)  $ are critical spaces  for the system $(S_1)$, in order to solve the problem of completeness of $\dot H^1$, the authors often use the inhomogeneous space $H^1$, it is also called critical space. In the periodic case ( see \cite{CCW}) the authors proved that if $\theta^0\in H^1\cap L^\infty$ such that $\|\theta^0\|_{L^\infty}<c_\infty\kappa$, then we get a global solution satisfying $\|\theta(t)\|_{H^1}\leq \|\theta^0\|_{H^1}$ for all $t\geq0$. Moreover, if in addition we have $\theta^0\in H^2$, they obtain an exponential decreasing property $\|\theta(t)\|_{H^2}\leq\|\theta^0\|_{H^2} e^{-c_\infty^{-1}(c_\infty\kappa-\|\theta^0\|_{L^\infty})t}$. In \cite{DC}, Dong a shows the global existence in $H^1(\mathbb T^2)$ with the only condition $\theta^0\in H^1$. The case of the whole plane $\R^2$ was treated by Dong in his good work \cite{DC}, more precisely he proved that there is global existence in $H^1(\R^2)$ if the initial condition is only small in space $\dot H^1(\R^2)$.\\

Our idea, inspired by \cite{CCW}, is to study the system $(S_1)$ in the critical space $X^0$ which is a little smaller than $L^\infty$.
In this paper we show some uniform regularity of the global solution: $\theta\in {\widetilde{L^\infty}}([0,T],X^0(\R^2))$ for any time $T$. This property helps us to show smallness at infinity and the analyticity of such a solution. Finally we prove the relationship between the  decreasing at infinity and the small frequencies of initial condition. To simplify the calculus and some steps of proofs of our results, we assume that $\kappa=1$, then we get
$$ (S_2)\begin{cases}
  \partial_t\theta+|D|\theta+u_\theta.\nabla\theta  = 0\;\;{\rm in}\;\mathbb{R}^{+*}\times\mathbb{R}^2\\
u_\theta=\mathcal{R}^\bot\theta=(-R_2\theta,R_1\theta)\\
\theta_{/_{t=0}}=\theta^0.
   \end{cases}$$
Our first result is the following.
\begin{theorem}\label{th11} Let $\theta^0\in X^0(\mathbb{R}^2)$ such that $\|\theta^0\|_{X^0}<1$, then there exists a unique solution $\theta$ of the system $(S_2)$ in the space $C_b([0,\infty),X^0(\mathbb{R}^2))$. Furthermore we have $\nabla \theta \in L^1([0,\infty),X^0(\mathbb{R}^2))$. Precisely
\begin{equation}\label{eqth1}\|\theta(t)\|_{ X^0}+(1-\|\theta^0\|_{X^0})\int_0^t\|\theta(z)\|_{X^1}dz\leq\|\theta^0\|_{X^0}.\end{equation}\end{theorem}
\begin{rem}\label{rem12} a) The local existence and uniqueness is due to the Fixed Point Theorem, but the global existence is due to the following inequality
\begin{equation}\label{remeq1}\|\theta\|_{X^0}+\int_0^t\|\theta\|_{ X^1}\leq\|\theta^0\|_{X^0}+\|\theta\|_{L^\infty([0,t],X^0)}\int_0^t\|\theta\|_{X^1},\end{equation}
and the smallness of the initial data.\\
b) If the global solution $\theta$ satisfies $\|\theta(t)\|_{ X^0}\leq \frac{1+\|\theta^0\|_{X^0}}{2}$ for all $t\geq0$, by using (\ref{remeq1}) we get
\begin{equation}\label{remeq2}\|\theta(t)\|_{ X^0}\leq \|\theta^0\|_{X^0},\;\forall t\geq 0.\end{equation}
Again by using equation (\ref{remeq1}) and (\ref{remeq2}) we get (\ref{eqth1}).\\
c) To show global existence we did not end up showing an explosion result, but we introduced a time noted $T_*$ which is not necessarily the maximum time of existence. The continuity of $(t\mapsto\|\theta(t)\|_{X^0})$ at 0 and the inequality (\ref{remeq1}) show that we were able to extend the solution beyond $T_*$, which gives the global solution.
\end{rem}
Our second result characterizes the regularizing effect of $(S_2)$ equations where we used the following notations
$$\mathfrak{X}^0=\{f\in C_b(\mathbb{R}^+,X^0(\R^2));\int_{\mathbb{R}^2}\Big[\sup_{0\leq t<\infty}|\hat{f}(t,\xi)|\Big]d\xi<\infty\}$$\;\;{\rm and}\;\;$$\mathfrak{X}^1=L^1(\mathbb{R}^+,X^1(\R^2)).$$
\begin{theorem}\label{th12}
  For any initial data in $X^0(\R^2)$ with $\|\theta^0\|_{X^0}<1/4$, there exists a unique global solution $\theta\in\mathfrak{X}^0\cap\mathfrak{X}^1$ such that $$\|\theta\|_{\mathfrak{X}^0}+\|\theta\|_{\mathfrak{X}^1}\leq4\|\theta^0\|_{X^0}.$$
\end{theorem}
By using the above results we can prove the decrease of global solution in large time towards 0. Precisely, we have the following result.
\begin{theorem}\label{th13} Let $\theta^0\in X^0(\mathbb{R}^2)$ such that $\|\theta^0\|_{X^0}<1$, then the global solution $\theta$ of the system $(S_2)$ given by Theorem \ref{th11} satisfies
\begin{equation}\label{eqcor1}\lim_{t\rightarrow+\infty}\|\theta(t)\|_{ X^0}=0.\end{equation}
\end{theorem}
Now, we state the second principle result of our work.
\begin{theorem}\label{th14} Let $\theta^0\in X^0(\R^2)$ such that $\|\theta^0\|_{\mathcal{X}^0}<1/8$, then the solution of $(S_2)$  by Theorem \ref{th11} satisfies $$\|e^{\frac{t}{2}|D|}\theta\|_{\mathfrak{X}^0}+\|e^{\frac{t}{2}|D|}\theta\|_{\mathfrak{X}^1}\leq12\|\theta^0\|_{{X}^0}.$$
\end{theorem}
\begin{rem} a) If $\|\theta^0\|_{X^0}<1$, by Theorem \ref{th13} and Theorem \ref{th14}, there is a time $t_0$ such that the unique global solution $\theta$ given by Theorem \ref{th11} satisfies
$$e^{\frac{t}{2}|D|}\theta(t_0+.)\in\mathfrak{X}^0\cap\mathfrak{X}^1,$$
and $$\|e^{\frac{t}{2}|D|}\theta(t_0+.)\|_{\mathfrak{X}^0}
+\|e^{\frac{t}{2}|D|}\theta(t_0+.)\|_{\mathfrak{X}^1}\leq12\|\theta(t_0)\|_{X^0}.$$
b) We could improve the index of analyticity in the following way: For $0<\sigma<1$, there is $C_\sigma>0$ such that
$$\|e^{\sigma t|D|}\theta\|_{\mathfrak{X}^0}+\|e^{\sigma t|D|}\theta\|_{\mathfrak{X}^1}\leq C_\sigma\|\theta^0\|_{{X}^0},$$
under the condition $\|\theta^0\|_{X^0}<\varepsilon_\sigma$, for some small positif $\varepsilon_\sigma$. The case $\sigma = 1$ is critical in the right way $|D|e^{-(1-\sigma)t|D|}\theta^0\in L^1([0,\infty),X^1(\R^2))$. Generally, this property is valuable on $[T,\infty)$ for some $T\geq0$. Indeed: By Theorem \ref{th13} there is a time $t_\sigma$ such that $\|\theta(t_\sigma)\|_{X^0}<\varepsilon_\sigma$, then $$\|e^{\sigma t|D|}\theta(t_\sigma+.)\|_{\mathfrak{X}^0}
+\|e^{\sigma t|D|}\theta(t_\sigma+.)\|_{\mathfrak{X}^1}\leq C_\sigma\|\theta(t_\sigma)\|_{X^0}.$$
\end{rem}
Our last result is the asymptotic study of the global solution given by Theorem \ref{th11} with additional conditions on the small frequencies of the initial condition.
\begin{theorem}\label{th15} Let $0\leq \delta<1$ and $\theta^0\in{{X}}^0(\mathbb{R}^2)\cap{{X}}^{-\delta}(\mathbb{R}^2)$ such that $\|\theta^0\|_{X^0}<2^{-(3-\delta)}$, then there exists a unique solution $\theta$ of the system $(QG)$ in the space $C_b([0,\infty),{{X}}^0(\mathbb{R}^2)\cap{{X}}^{-\delta}(\mathbb{R}^2)$. Furthermore we have $\nabla \theta \in L^1([0,\infty),{{X}}^0(\mathbb{R}^2)\cap{{X}}^{-\delta}(\mathbb{R}^2))$. Moreover
 \begin{equation}\label{eq1th5}\|\theta(t)\|_{X^0}=o(t^{-\delta}),\;t\rightarrow\infty.\end{equation}
\end{theorem}
\begin{rem} a) The local existence and uniqueness is due to the Fixed Point Theorem, but the global existence is due the following inequality
$$\|\theta\|_{X^{-\delta}}+\int_0^t\|\theta\|_{X^{1-\delta}}\leq\|\theta^0\|_{X^0}+2^{1-\delta}\|\theta\|_{L^\infty([0,t],X^0)}\int_0^t\|\theta\|_{X^{1-\delta}}.$$
b) The result (\ref{eq1th5}) prove relationship between the quality of decreasing of $\|\theta(t)\|_{X^0}$ and the small frequencies of $\theta^0$.
\end{rem}
The paper is organized in the following way: In section $2$, we give some notations and important preliminary results. Section $3$ we prove the global existence if $\|\theta^0\|_{X^0}<1$ . The proof used Fixed Point Theorem and standard Fourier techniques. In section $4$, we give a proof of Theorem \ref{th12}. Section 5 is devoted to prove the long time decay given by Theorem \ref{th13}. In section $6$, we prove the analyticity property of $\theta$, which helps us to prove the last result. Section 7 is devoted to prove long time decay of the global solution if the initial condition has some regularity with respect to small frequencies. Finally in section 8, we give some general remarks, among others we give some results on the periodic case.
\section{Notations and preliminaries results}
\subsection{Notations}
\begin{enumerate}
\item[$\bullet$] The Fourier transformation is normalized as
$$
\mathcal{F}(f)(\xi)=\widehat{f}(\xi)=\int_{\mathbb R^2}\exp(-ix.\xi)f(x)dx,\,\,\,\xi=(\xi_1,\xi_2)\in\mathbb R^2.
$$
\item[$\bullet$] The inverse Fourier formula is
$$
\mathcal{F}^{-1}(g)(x)=(2\pi)^{-2}\int_{\mathbb R^2}\exp(i\xi.x)g(\xi)d\xi,\,\,\,x=(x_1,x_2)\in\mathbb R^2.
$$
\item[$\bullet$] The convolution product of a suitable pair of function $f$ and $g$ on $\mathbb R^2$ is given by
$$
(f\ast g)(x):=\int_{\mathbb R^2}f(y)g(x-y)dy.
$$
\item[$\bullet$] If $f,g:\R^2\rightarrow\R^2$ are two functions, we set
$$
u_f.\nabla g={\rm div}\,(gu_f).
$$
\item[$\bullet$] Let $(B,||.||)$ be a functional Banach space, $1\leq p \leq\infty$ and  $T>0$. We define $L^p_T(B)$ the space of all
measurable functions $[0,t]\ni t\mapsto f(t) \in B$ such that $t\mapsto||f(t)||\in L^p([0,T])$.\\
\item [$\bullet$] $\|f\|_{\widetilde{L^\infty}(I,B)}=\|\,\|f(.)\|_{L^\infty(I)}\|_{B}$.\\
\item[$\bullet$] The Sobolev space $H^s(\mathbb R^2)=\{f\in \mathcal S'(\mathbb R^2);\;(1+|\xi|^2)^{s/2}\widehat{f}\in L^2(\mathbb R^2)\}$.\\
\item[$\bullet$] The homogeneous Sobolev space $\dot H^s(\mathbb R^2)=\{f\in \mathcal S'(\mathbb R^2);\;\widehat{f}\in L^1_{loc}\;{\rm and}\;|\xi|^s\widehat{f}\in L^2(\mathbb R^2)\}$.
\item[$\bullet$] The Fourier space $X^\sigma(\mathbb R^2)=\{f\in \mathcal S'(\mathbb R^2);\;\widehat{f}\in L^1_{loc}\;{\rm and}\;|\xi|^\sigma\widehat{f}\in L^1(\mathbb R^2)\}$. We recall that $X^{\sigma}(\mathbb{R}^2)$\; is a Banch space\;if \;and \;only\;if\;$\sigma\leq 0.$
\item[$\bullet$] For $f\in H^s(\R^2 )$ or $f\in X^\sigma(\R^2 )$, we have $|\mathcal F(u_f)(\xi)|=|\widehat{f}(\xi)|$, almost everywhere.
\item[$\bullet$] The periodic Fourier transformation is normalized as
$$
\mathcal{F}(f)(k)=\widehat{f}(k)=\int_{\mathbb T^2}\exp(-ik.x)f(x)dx,\,\,\,k=(k_1,k_2)\in\mathbb Z^2.
$$
\item[$\bullet$] The inverse Fourier formula in periodic case is
$$
\mathcal{F}^{-1}((a_k)_{k\in\Z^2})(x)=\sum_{k\in\Z}a_k e^{ik.x},\,\,\,x=(x_1,x_2)\in\mathbb T^2.
$$
\item[$\bullet$] We use often the following equalities: For $T>0$, we have
$$\{(z,t)\in[0,T]^2/\,0\leq z\leq t\}=\{(z,t)\in[0,T]^2/\,z\leq t\leq T\}$$
and
$$\{(z,z',t)\in[0,T]^3/\,0\leq z\leq t,\,\,0\leq z'\leq t\}=\{(z,z',t)\in[0,T]^3/\,\max(z,z')\leq t\leq T\}.$$

\end{enumerate}
\subsection{Preliminaries results}
We begin by recalling an inequality cited in \cite{JMa1} which helps us to prove the decrease of the global solution if the initial condition is more regular.
\begin{lem}(See \cite{JMa1})\label{lm 1.5}
 If $f:\mathbb{R}^+\rightarrow\mathbb{R}^+$ is a continuous function such that \\$f(t)\leq M_0+\varepsilon_1f(\varepsilon_2t)$, with $\varepsilon_1\in]0,1[$ and $\varepsilon_2\in[0,1]$, then $f$ is bounded and $\|f\|_{L^\infty(\R^+)}\leq M_0(1-\varepsilon_1)^{-1}.$
\end{lem}
\begin{lem}\label{lm 1.1}
 If $f,g \;\in X^0(\mathbb{R}^2)$, then $fg \;\in X^0(\mathbb{R}^2)$ and $$\|fg\|_{X^0}\leq\|f\|_{X^0}\|g\|_{X^0}.$$
 Furthermore, if $f,g \;\in X^0(\mathbb{R}^2)\cap X^1(\mathbb{R}^2)$, then
  $$\|fg\|_{X^1}\leq \|f\|_{X^0}\|g\|_{X^1}+\|f\|_{X^1}\|g\|_{X^0}.$$
\end{lem}
{\bf Proof.} We have\begin{eqnarray*}
        \|fg\|_{X^0}&=& \int_{\R^2}|\mathcal F(fg)(\xi)|d\xi \\
         &\leq& \int_{\mathbb{R}^2}|\widehat{f}*\widehat{g}(\xi)|d\xi \\
         &\leq& \|\widehat{f}\|_{L^1}\|\widehat{g}\|_{L^1} \\
         &\leq& \|f\|_{X^0}\|g\|_{X^0}.
      \end{eqnarray*}
      Using $|\xi|\leq |\xi-\eta|+|\eta|$, we get
\begin{eqnarray*}
        \|fg\|_{X^1}&=& \int_{\R^2}|\xi|.|\mathcal F(fg)(\xi)| \\
         &\leq& \int_{\mathbb{R}^2}|\xi|.|\widehat{f}*\widehat{g}(\xi)| \\
          &\leq& \int_{\mathbb{R}^2}|\xi|\int_{\R^2}|\widehat{f}(\xi-\eta)|.|\widehat{g}(\eta)| \\
         &\leq& \int_{\mathbb{R}^2}\int_{\R^2}|\xi-\eta|.|\widehat{f}(\xi-\eta)|.|\widehat{g}(\eta)|
         +\int_{\mathbb{R}^2}\int_{\R^2}|\widehat{f}(\xi-\eta)|.|\eta|.|\widehat{g}(\eta)| .
      \end{eqnarray*}
       Young's inequality gives
        $$\|fg\|_{X^1} \leq \|f\|_{X^1}\|g\|_{X^0}+\|f\|_{X^0}\|g\|_{X^1}.$$
\begin{lem}\label{lm 1.2}
  Let $\theta\in L_T^\infty(X^0(\mathbb{R}^2))\cap L_T^1(X^{1}(\mathbb{R}^2))$, then
  $$\|\int_{0}^{t}e^{-(t-z)|D|}(u_\theta\nabla\theta)dz\|_{X^0}\leq \|\theta\|_{L_T^\infty(X^0)}\|\theta\|_{L_T^1(X^1)}.$$
\end{lem}
{\bf Proof.}  We have $$\|\int_{0}^{t}e^{-(t-z)|D|}(u_\theta\nabla \theta)dz\|_{X^0}\leq \int_{0}^{t}\|u_\theta\nabla \theta\|_{X^0}dz.$$
  Using Lemma \ref{lm 1.1} and the fact that $\|u_\theta\|_{X^0}=\|\theta\|_{X^0}$ and $\hat{u_\theta}(\xi)=i(\frac{\xi_2}{|\xi|},-\frac{\xi_1}{|\xi|})\hat{\theta}$, we obtain\begin{eqnarray*}
      \|\int_{0}^{t}e^{-(t-z)|D|}(u_\theta\nabla \theta)dz\|_{X^0} &\leq &  \|u_\theta\|_{L_T^\infty(X^0)}\|\theta\|_{L_T^1(X^1)} \\
       &\leq&  \|\theta\|_{L_T^\infty(X^0)}\|\theta\|_{L_T^1(X^1)}.
                                       \end{eqnarray*}
\begin{lem}\label{lm 1.3}
  Let $\theta\in L_T^\infty(X^0(\mathbb{R}^2))\cap L_T^1(X^{1}(\mathbb{R}^2))$, then
  $$\int_{0}^{T}\|\int_{0}^{t}e^{-(t-z)|D|}(u_\theta\nabla \theta)dz\|_{X^1}dt\leq \|\theta\|_{L_T^\infty(X^0)}\|\theta\|_{L_T^1(X^1)}.$$
\end{lem}
{\bf Proof.} We have \begin{eqnarray*}
         \int_{0}^{T}\|\int_{0}^{t}e^{-(t-z)|D|}(u_\theta\nabla \theta)dz\|_{X^1}dt &\leq& \int_{0}^{T}\int_{0}^{t}\|e^{-(t-z)|D|}(u_\theta\nabla \theta)(z)\|_{X^1}dzdt \\
          &\leq& \int_{0}^{T}\int_{0}^{t}\int_{\mathbb{R}^2}e^{-(t-z)|\xi|}|\xi||\mathcal{F}(u_\theta\nabla \theta)(z,\xi)|d\xi dzdt \\
          &\leq& \int_{\mathbb{R}^2}|\xi|\Big(\int_{0}^{T}\int_{0}^{t}e^{-(t-z)|\xi|}|\mathcal{F}(u_\theta\nabla \theta)(z,\xi)| dzdt\Big)d\xi.
       \end{eqnarray*}
 Integrating the function $e^{-(t-z)|\xi|}$ twice with $z\in [0,t]$ and $t\in [0,T]$, we obtain
 $$\int_{0}^{T}\int_{0}^{t}e^{-(t-z)|\xi|}|\mathcal{F}(u_\theta\nabla\theta)(z,\xi)| dzdt=\int_{0}^{T}|\xi|\mathcal{F}(u_\theta\nabla \theta)(z,\xi)|\Big(\int_{z}^{T}e^{-(t-z)|\xi|}|dt\big)dz.$$
 Then
 \begin{eqnarray*}
   \int_{0}^{T}\|\int_{0}^{t}e^{-(t-z)|D|}(u_\theta\nabla \theta)dz\|_{X^1}dt &\leq& \int_{\mathbb{R}^2}|\xi|\Big(\int_{0}^{T}\Big[\int_{z}^{T}e^{-(t-z)|\xi|}dt\Big]|\mathcal{F}(u_\theta\nabla \theta)(z,\xi)| dz\Big)d\xi \\
    &\leq& \int_{\mathbb{R}^2}|\xi|\Big(\int_{0}^{T}\Big[\frac{1-e^{-(T-z)|\xi|}}{|\xi|}\Big]|\mathcal{F}(u_\theta\nabla \theta)(z,\xi)| dz\Big)d\xi \\
    &\leq& \int_{\mathbb{R}^2}\Big(\int_{0}^{T}|\mathcal{F}(u_\theta\nabla v)(z,\xi)| dz\Big)d\xi \\
    &\leq& \int_{0}^{T}\|u_\theta\nabla \theta\|_{X^0}dz.
 \end{eqnarray*}
 This implies
 \begin{eqnarray*}
   \int_{0}^{T}\|\int_{0}^{t}e^{(t-z)|D|}(u_\theta\nabla \theta)dz\|_{X^1}dt &\leq& \int_{0}^{T}\|u_\theta\nabla\theta\|_{X^0}\\
    &\leq& \|\theta\|_{L_T^\infty(X^0)}\|\theta\|_{L_T^1(X^1)}.
 \end{eqnarray*}
\begin{lem}\label{lem 5.2} For $\sigma>-1$, we have $L^2(\R^2)\cap X^{\sigma+1}(\R^2)\hookrightarrow X^\sigma(\R^2)$. Moreover, if $f\in L^2\cap X^{\sigma+1}$, then $$\|f\|_{X^\sigma}\leq(\sqrt{2\pi}+1)\|f\|_{L^2}^{\frac{1}{\sigma+2}}\|f\|_{X^{\sigma+1}}^{\frac{\sigma+1}{\sigma+2}}.$$
\end{lem}
{\bf Proof.}
For $\lambda>0$, we can write $$\|f\|_{X^\sigma}=\underbrace{\int_{|\xi|<\lambda}|\xi|^\sigma|\hat{f}(\xi)|}_{=I_\lambda}
+\underbrace{\int_{|\xi|>\lambda}|\xi|^\sigma|\hat{f}(\xi)|}_{=J_\lambda}.$$
We do some estimations, to get \begin{eqnarray*}
          I_\lambda &\leq& \int_{|\xi|<\lambda}|\xi|^\sigma|\hat{f}(\xi)| \\
           &\leq& (\int_{|\xi|<\lambda}|\xi|^{2\sigma} d\xi)^{1/2}(\int_{|\xi|<\lambda}|\hat{f}(\xi)|^2d\xi)^{1/2}  \\
           &\leq& \sqrt{2\pi}(\int_{0}^{\lambda}r^{2\sigma+1}dr)^{1/2}\|f\|_{L^2}\\
           &\leq& \sqrt{2\pi}\lambda^{\sigma+1}\|f\|_{L^2},
        \end{eqnarray*}
and \begin{eqnarray*}
      J_\lambda &=& \int_{|\xi|>\lambda}\frac{1}{|\xi|}|\xi||\hat{f}(\xi)| \\
       &\leq& \frac{1}{\lambda}\|f\|_{X^{\sigma+1}}.
    \end{eqnarray*}
Choosing $\lambda=\Big(\frac{\|f\|_{X^{\sigma+1}}}{\|f\|_{L^2}}\Big)^{1/{\sigma+2}}$, we obtain the desired result.

\begin{lem}\label{lm 1.40}
  Let $\sigma\geq0$ and $a,b\in X^0(\R^2)\cap X^\sigma(\R^2)$, then $ab\in X^\sigma(\R^2)$ and
$$\|ab\|_{X^\sigma}\leq 2^\sigma(\|a\|_{X^0}\|b\|_{X^\sigma}+\|a\|_{X^\sigma}\|b\|_{X^0}).$$
\end{lem}
{\bf Proof.} We have
 \begin{align*}
\|ab\|_{X^\sigma}&\leq \int_\xi  |\xi|^\sigma |\widehat{ab}(\xi)|d\xi\\
                        &\leq  \int_\xi \int_\eta |\xi|^{\sigma} |\widehat{a}(\xi-\eta)||\widehat{b}(\eta)|d\eta d\xi.
\end{align*}
And we have\begin{align*}
             |\xi| & \leq |\xi-\eta|+|\eta| \\
              & \leq 2\max(|\xi-\eta|,|\eta|).
           \end{align*}
 Then\begin{equation}\label{eq 1}
 |\xi|^{\sigma} \leq 2^{\sigma}|\xi-\eta|^{\sigma}+2^{\sigma}|\eta|^{\sigma},
     \end{equation}
     Similarly to the proof of lemma \ref{lm 1.1} we get
    $$\|ab\|_{X^\sigma}\leq  2^{\sigma}\int_{\mathbb{R}^2} \int_{\mathbb{R}^2} |\xi-\eta|^{\sigma} |\widehat{a}(\xi-\eta)||\widehat{b}(\eta)|d\eta d\xi+2^{\sigma}\int_{\mathbb{R}^2} \int_{\mathbb{R}^2}  |\widehat{a}(\xi-\eta)||\eta|^{\sigma}|\widehat{b}(\eta)|d\eta d\xi.$$
Which ends the proof.

\begin{lem}\label{lm 1.4}
  Let $\sigma\geq-1$, $a,b\in C_T(X^0\cap X^\sigma)\cap L^1_T(X^1\cap X^{\sigma+1})$ and $$Q(a,b)=\int^t_0e^{-(t-s)|D|}{\rm div}(bu_a).$$ Then
$$\begin{array}{lcl}
\|Q(a,b)\|_{L^\infty_T(X^\sigma)}&\leq&2^{\sigma+1}
\Big[\|a\|_{L^\infty_T(X^0)}\|b\|_{L^1_T(X^{\sigma+1})}+\|a\|_{L^1_T(X^{\sigma+1})}\|b\|_{L^\infty_T(X^0)}\Big]\\\\
\|Q(a,b)\|_{L^1_T(X^{\sigma+1})}&\leq&2^{\sigma+1}\Big[\|a\|_{L^\infty_T(X^0)}\|b\|_{L^1_T(X^{\sigma+1})}+\|a\|_{L^1_T(X^{\sigma+1})}\|b\|_{L^\infty_T(X^0)}\Big].\end{array}$$
Particularly
$$\begin{array}{lcl}
\|Q(a,a)\|_{L^\infty_T(X^\sigma)}&\leq&2^{\sigma+2}\|a\|_{L^\infty_T(X^0)}\|a\|_{L^1_T(X^{\sigma+1})}\\\\
\|Q(a,a)\|_{L^1_T(X^{\sigma+1})}&\leq&2^{\sigma+2}\|a\|_{L^\infty_T(X^0)}\|a\|_{L^1_T(X^{\sigma+1})}.\end{array}$$
\end{lem}
{\bf Proof.} We have
 \begin{align*}
\|Q(a,b)(t)\|_{X^\sigma}&\leq \int_0^t\int_\xi e^{-(t-s)|\xi|} |\xi|^\sigma |\widehat{{\rm div}(u_a\otimes b)}(s,\xi)|d\xi ds\\
                        &\leq  \int_0^t\int_\xi e^{-(t-s)|\xi|} |\xi|^{\sigma+1} |\widehat{(u_a\otimes b)}(s,\xi)|d\xi ds\\
                        &\leq  \int_0^t\int_\xi \int_\eta e^{-(t-s)|\xi|} |\xi|^{\sigma+1} |\widehat{u_a}(s,\xi-\eta)||\widehat{b}(s,\eta)|d\eta d\xi ds.
\end{align*}
By inequality (\ref{eq 1}), we have  \begin{align*}
\|Q(a,b)(t)\|_{X^\sigma}&\leq 2^{\sigma+1}\int^T_0 \int_\xi\int_\eta|\eta|^{\sigma+1}|\widehat{a}(s,\xi-\eta)||\widehat{b}(s,\eta)|+ |\xi-\eta|^{\sigma+1}|\widehat{a}(s,\xi-\eta)||\widehat{b}(s,\eta)|d\eta d\xi ds\\
                       &\leq 2^{\sigma+1}\big[\|a\|_{L^\infty_T(X^0)}\|b\|_{L^1_T(X^{\sigma+1})}+
                       \|a\|_{L^1_T(X^{\sigma+1})}\|b\|_{L^\infty_T(X^0)}\big].
 \end{align*}
Similarly, we have \begin{align*}
\|Q(a,b)(t)\|_{L^1_T(X^{\sigma+1})}&\leq \int_\xi\int_\eta\int^T_0\int^t_0 e^{-(t-s)|\xi|}|\xi|^{\sigma+2}|\widehat{a}(s,\xi-\eta)||\widehat{b}(s,\eta)|ds dt d\eta d\xi\\
&\leq \int_\xi\int_\eta\int^T_0\int^T_s e^{-(t-s)|\xi|}|\xi|^{\sigma+2}|\widehat{a}(s,\xi-\eta)||\widehat{b}(s,\eta)|dt ds d\eta d\xi\\
                                   &\leq \int_\xi \int_\eta \int^T_0 (1-e^{-(T-s)|\xi|}) |\xi|^{\sigma+1}|\widehat{a}(s,\xi-\eta)||\widehat{b}(s,\eta)|ds d\eta d\xi\\
                                   &\leq 2^{\sigma+1}[\|a\|_{L^\infty_T(X^0)}
                                   \|b\|_{L^1_T(X^{\sigma+1})}+\|a\|_{L^1_T(X^{\sigma+1})}\|b\|_{L^\infty_T(X^0)}].
 \end{align*}

\section{Proof of Theorem \ref{th11}}
In this section we prove the Theorem \ref{th11}. In the first part, we show the local existence and the uniqueness, and the second part is devoted to prove the global existence and the estimate (\ref{eqth1}).
\subsection{Uniqueness and local existence when $\|\theta^0\|_{X^0}<1$}
For $T>0$, we define the norm on the space $C_T(X^0)\cap L^1_T(X^1)$
$$N(\theta)=\max(\|\theta\|_{L^\infty_T(X^0)},\|\theta\|_{L^1_T(X^{1})}).$$
Let
$$r_0=\frac{1+\|\theta^0\|_{X^0}}{2}\in(\|\theta^0\|_{X^0},1),$$
$r_1$ is a small positif real number to fixed later, and the closed subset of $C_T(X^0)\cap L^1_T(X^1)$:
 $$B_T=\{\theta\in C_T(X^0)\cap L^1_T(X^1)/\|\theta\|_{L_T^\infty(X^0)}\leq r_0\;\;{\rm and}\;\|\theta\|_{L^1_T(X^1)}\leq r_1\}.$$
Now, Consider the following application
 \begin{align*}
   \psi: & B_T\longrightarrow  C_T(X^0)\cap L^1_T(X^1)\\
    & \theta\longmapsto e^{-t|D|}a-\int_{0}^{t}e^{-(t-z)|D|}u\nabla \theta.
 \end{align*}
$\bullet$ We start by looking for conditions on $r_1$ and $T$ such that $\psi(B_T)\subset B_T$:\\
Estimate of $\psi(\theta)$ in $C_T(X^0)$: We have $$\|e^{-t|D|}\theta^0\|_{L^\infty_T(X^0)}\leq\|\theta^0\|_{X^0},$$
and
\begin{eqnarray*}
  \|\int_{0}^{t}e^{-(t-z)|D|}u\nabla \theta\|_{L^\infty_T(X^0)} &\leq& \|\theta\|_{L^\infty_T(X^0)}\|\theta\|_{L^1_T(X^1)}\\
   &\leq&  r_0r_1,
\end{eqnarray*}
which imply
$$\|\psi(\theta)\|_{L^\infty_T(X^0)}\leq \|\theta^0\|_{X^0}+ r_0r_1.$$
Then the first sufficient condition is %
$$(C_1)\;\;\;\;\;\;\|\theta^0\|_{X^0}+ r_0r_1\leq r_0\Longleftrightarrow0<r_1\leq\frac{r_0-\|\theta^0\|_{X^0}}{r_0}=\frac{1-\|\theta^0\|_{X^0}}{1+\|\theta^0\|_{X^0}} .$$
Estimate of $\psi(\theta)$ in $L^1_T(X^1)$: We have \begin{eqnarray*}
      \|e^{-t|D|}\theta^0\|_{L^1_T(X^1)} &\leq& \int_{\R^2}(1-e^{T|\xi|})|\hat{\theta^0}(\xi)|d\xi\\
      \|\int_{0}^{t}e^{-(t-z)|D|}u\nabla \theta\|_{L^1_T(X^1)} &\leq& \|\theta\|_{L^\infty_T(X^0)}\|\theta\|_{L^1_T(X^1)}\leq r_0r_1.
    \end{eqnarray*}
Let $T=T(\theta^0,r_{1})>0$ such that
$$(C_2)\;\;\;\;\;\;\int_{\R^2}(1-e^{T|\xi|})|\hat{\theta^0}(\xi)|d\xi\leq(1-r_0)r_1.$$
Then  if the conditions $(C_1)$ and $(C_2)$ are  satisfied, we obtain $\psi(B_T)\subset B_T$.\\
$\bullet$ Search for additional conditions on $r_1$ and $T$ such as  $\psi$ is a contraction on $B_T$.\\
Let $\theta_1,\theta_2\in B_T$, then $$\varphi(\theta_1)-\varphi(\theta_2)=-\int^t_0 e^{-(t-s)|D|}u_{\theta_1-\theta_2} \nabla \theta_1-\int^t_0 e^{-(t-s)|D|}u_{\theta_2} \nabla (\theta_1-\theta_2).$$
Then, by lemma \ref{lm 1.2}, we get \begin{align*}
          \|\varphi(\theta_1)-\varphi(\theta_2)\|_{L^\infty_T(X^0)} & \leq \|\theta_1-\theta_2\|_{L^\infty_T(X^0)}\|\theta_1\|_{L^1_T(X^{1})}+ \|\theta_2\|_{L^\infty_T(X^0)}\|\theta_1-\theta_2\|_{L^1_T(X^{1})}\\
           & \leq (r_0+r_1)N_1(\theta_1-\theta_2),
        \end{align*}
        and
\begin{align*}
         \int_0^T \|\varphi(\theta_1)-\varphi(\theta_2)\|_{X^1} dt& \leq \|\theta_1-\theta_2\|_{L^\infty_T(X^0)}\|\theta_1\|_{L^1_T(X^{1})}+ \|\theta_2\|_{L^\infty_T(X^0)}\|\theta_1-\theta_2\|_{L^1_T(X^{1})}\\
           & \leq (r_0+r_1)N(\theta_1-\theta_2),
        \end{align*}
Then the third sufficient condition is
$$(C_3)\;\;\;\;\;\;r_0+r_1<1.$$
Therefore, if $(C_1)$, $(C_2)$ and $(C_3)$ are satisfied we get $\psi(B_T)\subset B_T$ and $\psi$ is contraction on $B_T$. Fixed Point Theorem gives the existence and uniqueness of solution of $(S_2)$ in $C_T(X^0)\cap L^1_T(X^1)$.
\subsection{Global existence} Let $\theta\in C_T(X^0)\cap L^1_T(X^1)$ be the solution of $(S_2)$ given by the first step. Then
\begin{equation}\label{eq1111}
  \partial_t \hat{\theta} + |\xi| \hat{\theta} + \mathcal{F}(u_\theta. \nabla \theta) = 0.
\end{equation}
We multiply this equation by $\overline{\hat{\theta}}$, we obtain
\begin{equation}\label{eq1112}
    \partial_t \hat{\theta} \overline{\hat{\theta}} + |\xi| {|\hat{\theta}|}^2 + \mathcal{F}(u_\theta. \nabla \theta)\overline{\hat{\theta}}= 0.
\end{equation}
By the equation (\ref{eq1111}) we have
$$
    \partial_t \overline{\hat{\theta}} +  |\xi| \overline{\hat{\theta}} + \overline{\mathcal{F}(u_\theta. \nabla \theta)} + \overline{\mathcal{F}(u_\theta. \nabla \theta)}= 0,
$$
We multiply this equation by $\hat{\delta}$ we obtain
\begin{equation}\label{eq1113}
     \partial_t \overline{\hat{\theta}} \hat{\theta} + |\xi| {|\hat{\theta}|}^2 + \overline{\mathcal{F}(u_\theta. \nabla \theta)} \hat{\theta}= 0.
\end{equation}
The sum of equations (\ref{eq1112}) and (\ref{eq1113}) gives
$$ \partial_t {|\hat{\theta}|}^2 + 2  |\xi| {|\hat{\theta}|}^2 + 2 Re(\mathcal{F}(u_\theta. \nabla \theta) \overline{\hat{\theta}})= 0$$
and
$$ \partial_t {|\hat{\theta}|}^2 + 2  |\xi| {|\hat{\theta}|}^2 \leq 2 |\mathcal{F}(u_\theta. \nabla \theta)|. |\overline{\hat{\theta}})|.$$
For $\varepsilon > 0$, we have
$$ \partial_t {|\hat{\theta}|}^2  = \partial_t ({|\hat{\theta}|}^2 + \varepsilon^2) = 2 \sqrt{{|\hat{\theta}|}^2 + \varepsilon^2 }\partial_t \sqrt{{|\hat{\theta}|}^2 + \varepsilon^2 }.$$
Then
$$2 \partial_t  \sqrt{{|\hat{\theta}|}^2 + \varepsilon^2 } + 2 |\xi| \frac{{|\hat{\theta}|}^2}{\sqrt{{|\hat{\theta}|}^2 + \varepsilon^2 }} \leq  2 |\mathcal{F}(u_\theta. \nabla \theta)| \frac{{|\hat{\theta}|}}{\sqrt{{|\hat{\theta}|}^2 + \varepsilon^2 }}.$$
Integrating with respect to time, we have
$$\sqrt{{|\hat{\theta}(t)|}^2 + \varepsilon^2 } + \int_0^t |\xi|  \frac{{|\hat{\theta}|}^2}{\sqrt{{|\hat{\theta}|}^2 + \varepsilon^2 }} \leq |\widehat{\theta^0}(\xi)|+
 \int_0^t |\mathcal{F}(u_\theta. \nabla \theta)|.$$
By tending $\varepsilon \longrightarrow 0$, we obtain
$$ |\hat{\theta}(t)| + \int_0^t |\xi| |\hat{\theta}| \leq |\widehat{\theta^0}(\xi)|+
 \int_0^t |\mathcal{F}(u_\theta. \nabla \theta)|.$$
Integrating with respect to $\xi$ we obtain
 $$\begin{array}{lcl}
\displaystyle{\|\theta(t)\|}_{X^0} +  \int_0^t  {\|\theta\|}_{X^{1}}&\leq&\displaystyle\|\theta^0\|_{X^0}+\int_0^t  {\| u_\theta. \nabla \theta\|}_{X^0}\\
 &\leq&\displaystyle\|\theta^0\|_{X^0}+\int_0^t  \|\theta\|_{X^0}\|\theta\|_{X^1}\\
  &\leq&\displaystyle\|\theta^0\|_{X^0}+r_0\int_0^t \|\theta\|_{X^1}\\
 \end{array}$$
Then
$${\|\theta(t)\|}_{X^0} +  (1-r_0)\int_0^t  {\|\theta\|}_{X^{1}}\leq\displaystyle\|\theta^0\|_{X^0},$$
which implies ${\|\theta(t)\|}_{X^0}\leq \|\theta^0\|_{X^0}$. Similarly, by integrating over $[s,t]\subset[0,T]$, we get
$$\|\theta(t)\|_{X^0}\leq \|\theta(s)\|_{X^0}.$$
Particularly $(t\longrightarrow\|\theta(t)\|_{X^0})$ is decreasing over $[0,T]$.
$$\forall t\in[0,T],\;\;{\|\theta(t)\|}_{X^0}\leq\|\theta^0\|_{X^0}.$$
Let $T_*\in(0,+\infty]$ the maximal time such that $\theta\in C([0,T_*),X^0)$ is a solution of $(S_2)$ and
$$\forall t\in[0,T_*),\;\;\|\theta(t)\|_{X^0}\leq r_0.$$
We want to prove $T_*=+\infty$. By using the same method of the first step and using the decay of the function $(t\rightarrow\|\theta(t)\|_{X^0})$, we get $\theta\in L^1_{loc}([0,T_*),X^1)$ and
$$\begin{array}{lcl}
\displaystyle{\|\theta(t)\|}_{X^0} +  \int_0^t  {\|\theta\|}_{X^{1}}&\leq&\displaystyle\|\theta^0\|_{X^0}+\int_0^t \|\theta\|_{X^0}{\|\theta\|}_{X^{1}}\\
&\leq&\displaystyle\|\theta^0\|_{X^0}+r_0\int_0^t {\|\theta\|}_{X^{1}}.
\end{array}$$
Then
\begin{equation}\label{eqff1}{\|\theta(t)\|}_{X^0} +  (1-r_0)\int_0^t  {\|\theta\|}_{X^{1}}\leq\displaystyle\|\theta^0\|_{X^0}.\end{equation}
Therefore $\theta\in L^1([0,T_*),X^1)$. We assume the opposite such as $T_*<\infty$. For $0\leq t<t'<T_*$, we have
$$\begin{array}{lcl}
\displaystyle{\|\theta(t')-\theta(t)\|}_{X^0}&\leq&\displaystyle\int_t^{t'}\|\theta\|_{X^{1}}+\int_t^{t'}\|u_\theta.\nabla\theta\|_{X^{0}}\\
&\leq&\displaystyle\int_t^{t'}  {\|\theta\|}_{X^{1}}+\int_t^{t'} \|\theta\|_{X^0}{\|\theta\|}_{X^{1}}\\
&\leq&\displaystyle\int_t^{t'}  {\|\theta\|}_{X^{1}}+r_0\int_t^{t'}{\|\theta\|}_{X^{1}}\\
&\leq&(1+r_0)\displaystyle\int_t^{t'}  {\|\theta\|}_{X^{1}}.
\end{array}$$
As $\displaystyle\int_t^{t'}  {\|\theta\|}_{X^{1}}\longrightarrow0$ if $t,t'\longrightarrow T_*$, then there is $\theta_*\in X^0(\R^2)$ such that $$\|\theta_*\|_{X^0}\leq r_0<1,$$ and
$$\lim_{t\rightarrow T_*}\|\theta(t)-\theta_*\|_{X^0}=0.$$
Then, by applying again the first step to the following system
$$\begin{cases}
  \partial_t\gamma+|D|\gamma+u_\gamma\nabla\gamma = 0\\
\gamma_{/_{t=0}}=\theta_*,
   \end{cases}$$
we get a time $T_0>0$ and a solution $\gamma\in C_{T_0}(X^0)\cap L^1_{T_0}(X^1)$ of this system satisfying
$$\|\gamma(t)\|_{X^0}\leq \|\theta_*\|_{X^0},\;\forall t\in [0,T_0].$$
The function $$f(t)=\left\{\begin{array}{l}
\theta(t)\;{\rm if}\;t\in[0,T_*)\\
\gamma(t-T_*)\;{\rm if}\;t\in[T_*,T_*+T_0]
\end{array}\right.$$ then defines an extension of the solution $\theta$ over $[0,T_*+T_0]$, and satisfying
$$\|f(t)\|_{X^0}\leq r_0,\;\forall t\in[0,T_*+T_0],$$
which contradicts our assumption, therefore $T_*=+\infty$. Finally equation (\ref{eqth1}) is given by Remark \ref{rem12} b) and (\ref{eqff1}).
\section{Proof of Theorem \ref{th12}}
This proof is inspired by that of the principle result of \cite{HB} and it's done in three steps. We begin by define the operator
$$\varphi:\mathfrak{X}^0\cap \mathfrak{X}^1\longrightarrow\mathfrak{X}^0\cap \mathfrak{X}^1.$$
$$\theta\longmapsto e^{-t|D|}\theta^0-\int_{0}^{t}e^{-(t-\tau)|D|}u_\theta.\nabla\theta(\tau)d\tau$$
{\bf {First step:}}  We first take the Fourier transformation to the above integral form of $\varphi(\theta)$, then we have
\begin{equation}\label{eq 2.1}
|\mathcal F(\varphi(\theta))(t,\xi)| \leq  e^{-t|\xi|}|\hat{\theta^0}(\xi)|+\int_{0}^{t}e^{-(t-\tau)|\xi|}|\hat{u}_\theta|\ast_\xi|\hat{\nabla\theta}|(\tau,\xi)d\tau \\
\end{equation}
then
\begin{equation}\label{eq 2.2}
  |\mathcal F(\varphi(\theta))(t,\xi)|\leq|\hat{\theta^0}(\xi)|+\Big[\sup_{0\leq z<\infty}|\widehat{\theta}(z,\xi)|\Big]\ast_\xi\Big[\int_{0}^{\infty}|\xi||\hat{\theta}(\tau,\xi)|d\tau\Big].
\end{equation}
By taking the $L^1$ norm in $\xi$ to (\ref{eq 2.2}) and applying Young inequality, we get
\begin{equation}\label{eq 2.3}
  \|\varphi(\theta)\|_{\mathfrak{X}^0}\leq\|\theta^0\|_{X^0}+\|\theta\|_{\mathfrak{X}^0}\|\theta\|_{\mathfrak{X}^1}.
\end{equation}
{\bf {Second step:}} Now we estimate $\theta$ in $\mathfrak{X}^1$. We multiply (\ref{eq 2.1}) by $|\xi|$, then
\begin{equation}\label{eq 2.4}
  |\xi|.|\mathcal F(\varphi(\theta))(t,\xi)|\leq  |\xi|e^{-t|\xi|}|\hat{\theta^0}(\xi)|+\int_{0}^{t}|\xi|e^{-(t-\tau)|\xi|}|\hat{u}\ast_\xi\hat{\nabla\theta}(\tau,\xi)|d\tau \\
\end{equation}
We take the $L^1$ norm in time to (\ref{eq 2.4}) and taking into account $\displaystyle\int_\tau^\infty|\xi|e^{-(t-\tau)|\xi|}dt=1$, and using Young's inequality we get
$$\begin{array}{lcl}
\displaystyle \int_{0}^{\infty}|\xi|.|\mathcal F(\varphi(\theta))(t,\xi)|dt&\leq& \displaystyle|\hat{\theta^0}(\xi)|+\int_{0}^{\infty}\int_0^t |\xi|e^{-(t-\tau)|\xi|}|\hat{u}_\theta|\ast_\xi\hat{\nabla\theta}|(\tau,\xi)d\tau dt\\
&\leq& \displaystyle|\hat{\theta^0}(\xi)|+\int_{0}^{\infty}\Big(\int_\tau^\infty |\xi|e^{-(t-\tau)|\xi|}dt\Big)|\hat{u}_\theta|\ast_\xi|\hat{\nabla\theta}|(\tau,\xi)d\tau\\
&\leq& \displaystyle|\hat{\theta^0}(\xi)|+\int_{0}^{\infty}|\hat{u}_\theta|\ast_\xi|\hat{\nabla\theta}|(\tau,\xi)d\tau,
\end{array}$$
which implies
\begin{equation}\label{eq 2.5}
 \int_{0}^{\infty}|\xi|.|\mathcal F(\varphi(\theta))(t,\xi)|dt  \leq |\hat{\theta^0}(\xi)|+\Big[\int_{0}^{\infty}|\xi||\hat{\theta}(\tau,\xi)|d\tau\Big]\ast\Big[\sup_{0\leq z<\infty}|\hat{\theta}(z,\xi)|\Big].
\end{equation}
We take the $L^1$ norm in $\xi$ to (\ref{eq 2.5}), we obtain
\begin{equation}\label{eq 2.6}
  \|\varphi(\theta)\|_{\mathfrak{X}^1}\leq\|\theta^0\|_{X^0}+\|\theta\|_{\mathfrak{X}^0}\|\theta\|_{\mathfrak{X}^1}.
\end{equation}
{\bf {Third step:}} Combining (\ref{eq 2.3}) and (\ref{eq 2.6}) we finally have
\begin{equation}\label{eq 2.7}
 \max(\|\varphi(\theta)\|_{\mathfrak{X}^0},\|\varphi(\theta)\|_{\mathfrak{X}^1})\leq\|\theta^0\|_{X^0}+ \max(\|\theta\|_{\mathfrak{X}^0},\|\theta\|_{\mathfrak{X}^1})^2.
\end{equation}
which implies, if $\|\theta^0\|_{X^0}<1/4$, we have $\varphi(P)\subset P$ where
$$P=\{f\in \mathfrak{X}^0\cap \mathfrak{X}^1/\, \max(\|f\|_{\mathfrak{X}^0},\|f\|_{\mathfrak{X}^1})\leq\frac{1-\sqrt{1-4\|\theta^0\|_{X^0}}}{2}=\frac{2\|\theta^0\|_{X^0}}{1+\sqrt{1-4\|\theta^0\|_{X^0}}}\}.$$
Finally the unique solution given by Theorem \ref{th11} is in $\mathfrak{X}^0\cap\mathfrak{X}^1$ for small initial data in $X^0$.
\section{Proof of Theorem \ref{th13}}
The idea of this proof is to write $\theta(t_0)$, for some $t_0$, as the sum of two functions, the first is small in norm $X^0$ and the second is more regular. For this, let $0<\varepsilon<1/8$. By inequality (\ref{eqth1}) there is a time $t_0\geq0$ such that
\begin{equation}\label{eqth13}\int_{0}^\infty\|\theta(t_0+z)\|_{X^1}dz=\int_{t_0}^\infty\|\theta(s)\|_{X^1}ds<\varepsilon/4.\end{equation}
Let $k\in \mathbb{N}$ be a large enough such that  $$\int_{A_k^c}|\hat{\theta}(t_0,\xi)|d\xi <\varepsilon/4,$$
where $$A_k=\{\xi\in\mathbb{R}^2/\;|\xi|\leq k\;{\rm and}\;|\hat{\theta}(t_0,\xi)|\leq k\}.$$
Now, put the following functions
$$\begin{array}{lcl}
a^0&=&\mathcal F^{-1}({\bf 1}_{A_k}(\xi)\widehat{\theta}(t_0,\xi))\\\\
b^0&=&\mathcal F^{-1}({\bf 1}_{A_k^c}(\xi)\widehat{\theta}(t_0,\xi))=\theta(t_0)-a^0.
\end{array}$$
Clearly, we have
$$\left\{\begin{array}{l}
a^0\in X^0(\R^2)\cap L^2(\R^2),\\\\
\|b^0\|_{X^0}<\varepsilon/4.
\end{array}\right.$$
Also, let $b\in C_b(\mathbb{R}^+,X^0)\cap L^1(\mathbb{R}^+,X^1)$ the unique solution of the following system
$$\begin{cases}
  \partial_tb+|D|b+u_b\nabla b  = 0\\
b_{/_{t=0}}=b^0,
   \end{cases}$$
given by Theorem\ref{th11} and Theorem\ref{th12}. Moreover, $b$ satisfies
 \begin{equation}\label{eqpb1}\|b(t)\|_{X^0}+(1-\|b^0\|_{X^0})\int_{0}^{t}\|b(z)\|_{X^1}dz\leq\|b^0\|_{X^0},\;\forall t\geq0,\end{equation}
 and
 \begin{equation}\label{eqpb2}\|b\|_{\mathfrak{X}^0}+\|b\|_{\mathfrak{X}^1}\leq4\|b^0\|_{X^0}.\end{equation}
By Theorem \ref{th11}, the function $(t\rightarrow \theta(t_0+t))$ is the unique solution of the following system
$$\begin{cases}
  \partial_t\gamma+|D|\gamma+u_\gamma\nabla\gamma  = 0,\\
\gamma_{/_{t=0}}=\theta(t_0).
   \end{cases}$$
We denote then $a(t):=\theta(t_0+t)-b(t)$ which implies $\theta(t_0+t)=a(t)+b(t)$. As $(t\rightarrow \theta(t_0+t))$ and $b$ are two elements of $C_b(\mathbb{R}^+,X^0)\cap L^1(\mathbb{R}^+,X^1)$, then $a\in C_b(\mathbb{R}^+,X^0)\cap L^1(\mathbb{R}^+,X^1)$ and we have
\begin{equation}\label{eqpa1}
\|a(t)\|_{X^0}\leq\|\theta(t_0+t)\|_{X^0}+\|b(t)\|_{X^0}\leq \|\theta^0\|_{X^0}+\|b^0\|_{X^0},\end{equation}
\begin{equation}\label{eqpa2}\int_{0}^{t}\|a\|_{X^1}\leq\displaystyle\int_{0}^{t}\|\theta(t_0+z)\|_{X^1}dz+\int_{0}^{t}\|b(z)\|_{X^1}dz\leq \displaystyle\varepsilon/2.\end{equation}
Moreover, $a$ is the unique solution of the following system
$$(S)\begin{cases}
    \partial_t a +|D|a +u_a\nabla a+u_a\nabla b+u_b\nabla a& =0 \\
    a_{/_{t=0}} =a^0.
  \end{cases}$$
\subsection{Estimate of $\|\widehat{a}\|_{\widetilde{L^\infty}([0,T],L^2(\R^2))}$}  The integral form of the system $(S)$
  $$a(t)=e^{-t|D|}a^0-\int_0^te^{-(t-z)|D|}u_a.\nabla (a+b)-\int_0^te^{-(t-z)|D|}u_b.\nabla a$$
implies, for $T>0$,
  $$\|\widehat{a}\|_{\widetilde{L^\infty}([0,T],L^2(\R^2)}\leq \|\widehat{a^0}\|_{L^2(\R^2)}+J_1+J_2,$$
  where
$$\begin{array}{lcl}
  J_1&=&\displaystyle\|\sup_{0\leq t\leq T}|\mathcal F\Big(\int_0^te^{-(t-z)|D|}u_a.\nabla (a+b)\Big)|\|_{L^2(\R^2)},\\
  J_2&=&\displaystyle\|\sup_{0\leq t\leq T}|\mathcal F\Big(\int_0^te^{-(t-z)|D|}u_b.\nabla a\Big)|\|_{L^2(\R^2)}.
\end{array}$$
Young's inequality implies
$$\begin{array}{lcl}
  J_1&\leq&\displaystyle\int_0^T\|\widehat{u_a}\ast_\xi\mathcal F(\nabla (a+b))\|_{L^2(\R^2)},\\
   &\leq&\displaystyle\int_0^T\|\widehat{u_a}(z)\|_{L^2(\R^2)}\|\mathcal F(\nabla (a+b))(z)\|_{L^1(\R^2)},\\
   &\leq&\displaystyle\int_0^T\|\widehat{a}(z)\|_{L^2(\R^2)}\|\theta(t_0+z)\|_{X^1}dz,\\
  &\leq&\displaystyle\|\widehat{a}\|_{L^\infty([0,T],L^2(\R^2))}\int_0^T\|\theta(t_0+z)\|_{X^1}dz\\
  &\leq&\displaystyle\frac{\varepsilon}{4}\|\widehat{a}\|_{L^\infty([0,T],L^2(\R^2))},\;(by\;(\ref{eqth13})).
\end{array}$$
In other hand, to estimate $J_2$, we write
$$\begin{array}{lcl}
\displaystyle|\int_0^te^{-(t-z)|\xi|}\mathcal F(u_b.\nabla a)(z,\xi)dz|&=&\displaystyle|\int_0^te^{-(t-z)|\xi|}\mathcal F({\rm div}\,(au_b)(z,\xi)dz|\\
&\leq&\displaystyle\int_0^te^{-(t-z)|\xi|}|\xi|.|\widehat{a}|\ast_\xi|\widehat{u_b}|(z,\xi)|dz\\
  &\leq&\displaystyle\int_0^te^{-(t-z)|\xi|}|\xi|\int_{\R^2}|\widehat{a}(z,\xi-\eta)|.|\widehat{b}|(z,\eta)|d\eta dz\\
  &\leq&\displaystyle\int_0^te^{-(t-z)|\xi|}|\xi|dz\int_{\R^2}\|\widehat{a}(\xi-\eta)\|_{L^\infty([0,t])}\|\widehat{b}(\eta)\|_{L^\infty([0,t])}d\eta \\
   &\leq&\displaystyle(1-e^{-t|\xi|})\int_{\R^2}\|\widehat{a}(\xi-\eta)\|_{L^\infty([0,T])}\|\widehat{b}(\eta)\|_{L^\infty([0,T])}d\eta \\
   &\leq&\displaystyle\int_{\R^2}\|\widehat{a}(\xi-\eta)\|_{L^\infty([0,T])}\|\widehat{b}(\eta)\|_{L^\infty([0,T])}d\eta \\
   &\leq&\displaystyle\|\widehat{a}(.)\|_{L^\infty([0,T])}\ast_\xi\|\widehat{b}(.)\|_{L^\infty([0,T])}(\xi).
\end{array}$$
Young inequality gives
$$\begin{array}{lcl}
\displaystyle\|\sup_{0\leq t\leq T}|\int_0^te^{-(t-z)|\xi|}\mathcal F(u_b.\nabla a)(z,\xi)dz|\|_{L^2}&\leq&\displaystyle\|\|\widehat{a}(.)\|_{L^\infty([0,T])}\ast_\xi\|\widehat{b}(.)\|_{L^\infty([0,T])}\|_{L^2}\\
&\leq&\displaystyle\|\widehat{a}\|_{\widetilde{L^\infty}([0,T],L^2)}\|\widehat{b}\|_{\widetilde{L^\infty}([0,T],L^1(\R^2))}\\
&\leq&\displaystyle\|\widehat{a}\|_{\widetilde{L^\infty}([0,T],L^2)}\|b\|_{\mathfrak{X}^0}\\
&\leq&\displaystyle4\|b^0\|_{X^0}\|\widehat{a}\|_{\widetilde{L^\infty}([0,T],L^2)},\;(by\;(\ref{eqpb2}))\\
&\leq&\displaystyle\varepsilon\|\widehat{a}\|_{\widetilde{L^\infty}([0,T],L^2)}.
\end{array}$$
Combining the above inequalities and the fact $0<\varepsilon<1/8$, we get $$\|\widehat{a}\|_{\widetilde{L^\infty}([0,T],L^2(\R^2))}\leq 2\|\widehat{a^0}\|_{L^2(\R^2)},$$ which implies
\begin{equation}\label{eql22}
\|\widehat{a}(t)\|_{L^2(\R^2)}\leq 2\|\widehat{a^0}\|_{L^2(\R^2)},\;\forall t\geq0.
\end{equation}
\subsection{Last step of the proof}
Now, we return to the proof of Theorem \ref{th13}. Let $$E_\varepsilon=\{t\geq 0/\|a(t)\|_{X^0}\geq\frac{\varepsilon}{2}\}.$$
Let $t\in E_\varepsilon$ so $\|a(t)\|_{X^0}\geq \frac{\varepsilon}{2}.$
Then, by Lemma \ref{lem 5.2} we have $$(\sqrt{2\pi}+1)\|a(t)\|_{L^2}^{1/2}\|a(t)\|_{X^1}^{1/2}\geq\frac{\varepsilon}{2}.$$
By squaring the last inequality and using inequality (\ref{eql22}), we obtain $$2(\sqrt{2\pi}+1)^2\|a^0\|_{L^2}\|a(t)\|_{X^1}\geq(\frac{\varepsilon}{2})^2.$$
Integrating on $E_\varepsilon$ we get $$2(\sqrt{2\pi}+1)^2\|a^0\|_{L^2}\int_{E_\varepsilon}\|a(t)\|_{X^1}\geq(\frac{\varepsilon}{2})^2\lambda_1(E_\varepsilon).$$
Consequently we obtain $$2(\sqrt{2\pi}+1)^2\|a^0\|_{L^2}\|a^0\|_{X^0}\geq(\frac{\varepsilon}{2})^2\lambda_1(E_\varepsilon).$$
Furthermore $$\lambda_1(E_\varepsilon)\leq \frac{2(\sqrt{2\pi}+1)^2\|a^0\|_{L^2}\|a^0\|_{X^0}}{(\frac{\varepsilon}{2})^2}=T_\varepsilon.$$
As $\lambda_1([0,T_\varepsilon+1])=T_\varepsilon+1>\lambda_1(E_\varepsilon)$, there exists $t_1\in[0,T_\varepsilon+1]\backslash E_\varepsilon$ such that
$\|a(t_1)\|_{X^0}< \frac{\varepsilon}{2}.$
Finally, as $\theta(t_0+t_1)=a(t_1)+b(t_1)$, we get \begin{align*}
          \|\theta(t_0+t_1)\|_{X^0} &= \|a(t_1)+b(t_1)\|_{X^0}\\
           & \leq \|a(t_1)\|_{X^0}+\|b(t_1)\|_{X^0}\\
           &  < \frac{\varepsilon}{2}+\|b^0\|_{X^0}\\
           &  < \frac{\varepsilon}{2}+\frac{\varepsilon}{2}=\varepsilon.
        \end{align*}
Now consider the system starting in $\theta(t_0+t_1)$ $$\begin{cases}
                                            \partial\Gamma+|D|\Gamma+u_\Gamma\nabla\Gamma=0 \\
                                            \Gamma_{/_{t=0}}=\theta(t_0+t_1).
                                          \end{cases}$$
By the uniqueness given by Theorem \ref{th11}, we get $\Gamma(s)=\theta(s+t_0)$ and $$\|\Gamma(s)\|_{X^0}\leq \|\Gamma(0)\|_{X^0} \;\forall \;s\geq 0.$$
Moreover, we obtain $$\|\theta(t)\|_{X^0} \leq\|\theta(t_0+t_1)\|_{X^0} <\varepsilon,\;\forall\; t\geq t_0+t_1,$$
which ends the proof.
\section{Proof of Theorem \ref{th14}} In this section we prove the analyticity of the solution of $(S_2)$ given by Theorem \ref{th11} if the initial data is small enough. We define the subspace of $\mathfrak{X}^0\cap \mathfrak{X}^1$
$$\mathfrak{Y}=\{f\in\mathfrak{X}^0\cap \mathfrak{X}^1/\;e^{\frac{t}{2}|D|}f(t) \in\mathfrak{X}^0\cap \mathfrak{X}^1\}.$$
Then, for $\theta\in \mathfrak{Y}$, we have
  $$\begin{array}{lcl}
 \displaystyle e^{\frac{t}{2}|\xi|}\mathcal F(\varphi({\theta}))(t,\xi)&=&\displaystyle e^{-\frac{t}{2}|\xi|}\hat{\theta^0}(\xi)-\int_{0}^{t}e^{-\frac{t}{2}|\xi|+\tau|\xi|}\hat{u}_\theta\ast_\xi\hat{\nabla\theta}(\tau,\xi)\\
  &=&\displaystyle e^{-\frac{t}{2}|\xi|}\hat{\theta^0}(\xi)-\int_{0}^{t}e^{-\frac{(t-\tau)}{2}|\xi|}e^{\frac{\tau}{2}|\xi|}\hat{u}_\theta\ast_\xi\hat{\nabla\theta}(\tau,\xi)\\
  &=&\displaystyle e^{-\frac{t}{2}|\xi|}\hat{\theta^0}(\xi)-\int_{0}^{t}e^{-\frac{(t-\tau)}{2}|\xi|}\int e^{\frac{\tau}{2}(|\xi|-|\eta|-|\xi-\eta|)}e^{\frac{\tau}{2}|\eta|}\hat{u}_{\theta}(\tau,\eta)e^{\frac{\tau}{2}|\xi-\eta|}\hat{\nabla\theta}(\tau,\xi-\eta)\\
  \end{array}$$
Since $e^{\frac{\tau}{2}(|\xi|-|\xi-\eta|-|\eta|)}\leq 1$ is uniformly bounded independently of $\tau,\,\xi$ and $\eta$, then
 $$|e^{\frac{t}{2}|\xi|}\mathcal F(\varphi({\theta}))(t,\xi)|\leq e^{-\frac{t}{2}|\xi|}\hat{\theta^0}(\xi)+\int_{0}^{t}e^{-\frac{(t-\tau)}{2}|\xi|}(e^{\frac{\tau}{2}|.|}|\hat{u}_\theta|)\ast_\xi(e^{\frac{\tau}{2}|.|}|\hat{\nabla\theta}|)(\tau,\xi)d\tau.$$
 Then we follow the proof of Theorem \ref{th12}, we get
 $$\|e^{\frac{t}{2}|D|}\mathcal F(\varphi({\theta}))\|_{\mathfrak{X}^0}\leq\|\theta^0\|_{X^0}+\|e^{\frac{t}{2}|D|}\theta\|_{\mathfrak{X}^0}\|e^{\frac{t}{2}|D|}\theta\|_{\mathfrak{X}^1},$$
 and
 $$\|e^{\frac{t}{2}|D|}\mathcal F(\varphi({\theta}))\|_{\mathfrak{X}^1}\leq2\|\theta^0\|_{X^0}+2\|e^{\frac{t}{2}|D|}\theta\|_{\mathfrak{X}^0}\|e^{\frac{t}{2}|D|}\theta\|_{\mathfrak{X}^1},$$
 which implies,
 $$\max\Big(\|e^{\frac{t}{2}|D|}\mathcal F(\varphi({\theta}))\|_{\mathfrak{X}^0},\frac{\|e^{\frac{t}{2}|D|}\mathcal F(\varphi({\theta}))\|_{\mathfrak{X}^1}}{2}\Big)\leq\|\theta^0\|_{X^0}+2\max\Big(\|e^{\frac{t}{2}|D|}\theta\|_{\mathfrak{X}^0},\frac{\|e^{\frac{t}{2}|D|}\theta\|_{\mathfrak{X}^1}}{2}\Big)^2.$$
Then, if $\|\theta^0\|_{X^0}<1/8$, we have $\varphi(Q)\subset Q$ where
$$Q=\{f\in \mathfrak{Y}/\, \max\Big(\|e^{\frac{t}{2}|D|}f\|_{\mathfrak{X}^0},\frac{\|e^{\frac{t}{2}|D|} f\|_{\mathfrak{X}^1}}{2}\Big)\leq\frac{1-\sqrt{1-8\|\theta^0\|_{X^0}}}{2}=\frac{4\|\theta^0\|_{X^0}}{1+\sqrt{1-8\|\theta^0\|_{X^0}}}\}.$$
Finally the unique solution given by Theorem \ref{th11} is in $\mathfrak{Y}$ and $\|e^{\frac{t}{2}|D|}\theta\|_{\mathfrak{X}^0}+\|e^{\frac{t}{2}|D|} \theta\|_{\mathfrak{X}^1}\leq12\|\theta^0\|_{X^0}$.
\section{Proof of Theorem \ref{th15}}
In this section we prove Theorem \ref{th14} where we assume that $\theta^0\in X^0(\R^2)\cap X^{-\delta}(\R^2)$, and
$$\|\theta^0\|_{X^0}<2^{-(3-\delta)}.$$
This proof is done in four steps
$$\begin{array}{l}
\bullet\;{\rm Uniqueness\;and\,local\;existence\;in\;}C_T(X^0\cap X^{-\delta})\cap L^1_T(X^1\cap X^{1-\delta})\\\\
\bullet\;{\rm Global\;existence\;in\;}C_b(\R^+,X^0\cap X^{-\delta})\cap L^1(\R^+,X^1\cap X^{1-\delta})\\\\
\bullet\;\lim_{t\rightarrow\infty}\|\theta(t)\|_{X^{-\delta}}=0\\\\
\bullet\;\lim_{t\rightarrow\infty}t^\delta\|\theta(t)\|_{X^0}=0.
\end{array}$$
\subsection{Step1: Uniqueness and local existence.} For $T>0$, we begin by equipping the space $C_T(X^0\cap X^{-\delta})\cap L^1_T(X^1\cap X^{1-\delta})$ of the norm
$$N'(f)=\max\Big(\|\theta\|_{L^\infty_T(X^0)},\|\theta\|_{L^\infty_T(X^{-\delta})},\|\theta\|_{L^1_T(X^1)},\|\theta\|_{L^1_T(X^{1-\delta})}\Big).$$
Now, consider the following subset of $C_T(X^0\cap X^{-\delta})\cap L^1_T(X^1\cap X^{1-\delta})$:
$$A= \{\theta\in C_T(X^0\cap X^{-\delta})\cap L^1_T(X^1\cap X^{1-\delta})/\;\left\{\begin{array}{l}
\|\theta\|_{L^\infty_T(X^0)}\leq r_0\\
\|\theta\|_{L^\infty_T(X^{-\delta})}\leq 2\|\theta^0\|_{X^{-\delta}}\\
\|\theta\|_{L^1_T(X^1)}\leq r_1\\
\|\theta\|_{L^1_T(X^{1-\delta})}\leq r_2.
\end{array}\right.\}.$$
$\bullet$ Firstly and for a good choice of $T,\,r_0,\,r_1$ and $r_2$, we want to show that $\varphi(A)\subset A$. We have \begin{align*}
\|\varphi(\theta)\|_{L^\infty_T(X^0)}&\leq \|\theta^0\|_{X^0}+\|\theta\|_{X^0}\|\theta\|_{L^1_T(X^{1})}\\
                               &\leq \|\theta^0\|_{X^0}+r_0r_1.
\end{align*}
Then if $r_0$ and $r_1$ satisfy $$0<r_1\leq \frac{r_0-\|\theta^0\|_{X^0}}{r_0},$$ we get $\|\varphi(\theta)\|_{L^\infty_T(X^0)}\leq r_0$.
Also we have \begin{align*}
\|\varphi(\theta)\|_{L^1_T(X^{1})}&\leq \int_{\mathbb{R}^2}(1-e^{-T|\xi|})|\hat{\theta^0}|d\xi+\|\theta\|_{L^\infty_T(X^0)}\|\theta\|_{L^1_T(X^{1})}\\
                            &\leq \int_{\mathbb{R}^2}(1-e^{-T|\xi|})|\hat{\theta^0}|d\xi+r_0r_1.
\end{align*}
Then, by Dominated Convergence Theorem if we choose $T$ such that $$\int_{\mathbb{R}^2}(1-e^{-T|\xi|})|\hat{\theta^0}|d\xi\leq r_1(1-r_0),$$
we get $\|\varphi(\theta)\|_{L^1_T(X^{1})}\leq r_1$. In other hand we have \begin{align*}
\|\varphi(\theta)\|_{L^\infty_T(X^{-\delta})}&\leq \|\theta^0\|_{X^{-\delta}}+2^{2-\delta}\|\theta\|_{X^0}\|\theta\|_{L^1_T(X^{1-\delta})}\\
                               &\leq \|\theta^0\|_{X^{-\delta}}+2^{2-\delta}r_0r_2.
\end{align*}
Then if $r_0$ and $r_2$ satisfies $$r_2\leq\frac{\|\theta^0\|_{X^{-\delta}}}{2^{2-\delta}r_0},$$
we get $\|\varphi(\theta)\|_{L^\infty_T(X^{-\delta})}\leq 2\|\theta^0\|_{X^{-\delta}}.$
Finally \begin{align*}
\|\varphi(\theta)\|_{L^1_T(X^{1-\delta})}&\leq \int_{\mathbb{R}^2}(1-e^{-T|\xi|})|\hat{\theta^0}||\xi|^{-\delta} d\xi+2^{2-\delta}\|\theta^0\|_{X^0}\|\theta\|_{L^1_T(X^{1-\delta})}\\
                               &\leq \int_{\mathbb{R}^2}(1-e^{-T|\xi|})|\hat{\theta^0}||\xi|^{-\delta} d\xi+2^{2-\delta}r_0r_2.
\end{align*}
Then if we choose $r_0\in(\|\theta^0\|_{X^0},2^{-(2-\delta)})$ and $T>0$ such that
$$\int_{\mathbb{R}^2}(1-e^{-T|\xi|})|\hat{\theta^0}||\xi|^{-\delta} d\xi\leq r_2(1-2^{2-\delta}r_0),$$
we get $\|\varphi(\theta)\|_{L^1_T(X^{\sigma+1})}\leq r_2$.\\\\
So if the above conditions are all satisfied, then $\varphi(A)\subset A$.\\\\
$\bullet$ Secondly we look for additional conditions on $r_0,\,r_1,\,r_2$ and $T$  so that $\varphi$ is a contraction on $A$.  For this let  $\theta_1,\theta_2\in A$, then $$\varphi(\theta_1)-\varphi(\theta_2)=-\int^t_0 e^{-(t-s)|D|}u_{\theta_1-\theta_2} \nabla \theta_1-\int^t_0 e^{-(t-s)|D|}u_{\theta_2} \nabla (\theta_1-\theta_2).$$
We have \begin{align*}
          \|\varphi(\theta_1)-\varphi(\theta_2)\|_{X^0} & \leq \|\theta_1-\theta_2\|_{L^\infty_T(X^0)}\|\theta_1\|_{L^1_T(X^{1})}+ \|\theta_2\|_{L^\infty_T(X^0)}\|\theta_1-\theta_2\|_{L^1_T(X^{1})}\\
           & \leq (r_0+r_1)N'(\theta_1-\theta_2).
        \end{align*}
We have \begin{align*}
          \|\varphi(\theta_1)-\varphi(\theta_2)\|_{X^{-\delta}} & \leq 2^{2-\delta} \|\theta_1-\theta_2\|_{L^\infty_T(X^0)}\|\theta_1\|_{L^1_T(X^{1-\delta})}+ 2^{2-\delta} \|\theta_2\|_{L^\infty_T(X^0)}\|\theta_1-\theta_2\|_{L^1_T(X^{1-\delta})}\\
           & \leq 2^{2-\delta}(r_0+r_2)N'(\theta_1-\theta_2),
        \end{align*}
Similarly, we obtain $$\|\varphi(\theta_1)-\varphi(\theta_2)\|_{L^1_T(X^{1})})\leq (r_0+r_1)N'(\theta_1-\theta_2),$$
and$$\|\varphi(\theta_1)-\varphi(\theta_2)\|_{L^1_T(X^{1-\delta})}\leq 2^{2-\delta}(r_0+r_2)N'(\theta_1-\theta_2).$$
Therefore $$N'(\varphi(\theta_1)-\varphi(\theta_2))\leq \max\Big(r_0+r_1,2^{2-\delta}(r_0+r_2)\Big)N'(\theta_1-\theta_2).$$
Then, $\varphi$ is a contraction on $A$ if we add the additional condition $$\max\Big(r_0+r_1,2^{2-\delta}(r_0+r_2)\Big)<1.$$
We can finish by applying the Fixed Point Theorem.
\subsection{Step2: Global existence in $X^0\cap X^{-\delta}$} Clearly by Theorem \ref{th11}, we get $\theta\in C_b(\R^+,X^0)\cap L^1(\R^+,X^1)$ and $\|\theta(t)\|_{X^0}\leq \|\theta^0\|_{X^0}$ for all time $t$. It remains to prove that $\theta\in C_b(\R^+,X^{-\delta})\cap L^1(\R^+,X^{1-\delta})$.
Now, define a time $$t^*=\sup\{T\geq0/\;\theta\in C([0,T],X^{-\delta}(\R^2))\cap L^1([0,T],X^{1-\delta})\}.$$ By the local existence step $t^*$ is well defined and $t^*\in(0,\infty]$. We want to prove that $t^*=\infty$. For this, suppose that $t^*<\infty$, then for $0\leq t<t^*$ we have
$$\|\theta(t)\|_{X^{-\delta}}+\int_{0}^{t}\|\theta\|_{X^{1-\delta}}\leq\|\theta^0\|_{X^{-\delta}}+2^{1-\delta}\|\theta^0\|_{X^0}\int_{0}^{t}\|\theta\|_{X^{1-\delta}}$$
which implies
\begin{equation}\label{eqt*}\|\theta(t)\|_{X^{-\delta}}+(1-2^{1-\delta}\|\theta^0\|_{X^0})\int_{0}^{t}\|\theta\|_{X^{1-\delta}}
\leq\|\theta^0\|_{X^{-\delta}}.\end{equation}
Then $\theta\in C([0,t^*),X^0)\cap L^\infty([0,t^*),X^{1-\delta})\cap L^1([0,t^*),X^{1-\delta})$. Now, prove that $\lim_{t\nearrow t^*}\theta(t)$ exists in $X^{1-\delta}.$ For $0<t<t'<t_*$, we have
$$\theta(t')-\theta(t)= -\int_{t}^{t'}|D|\theta-\int_t^{t'}u_\theta.\nabla \theta$$
and
$$\|\theta(t)-\theta(t')\|_{X^{-\delta}}\leq \int_{t}^{t'}\|\theta\|_{X^{1-\delta}}+2^{1-\delta}\int_{t}^{t'}\|\theta\|_{X^0}\|\theta\|_{X^{1-\delta}}\leq (1+2^{1-\delta}\|\theta^0\|_{X^0})\int_{t}^{t'}\|\theta\|_{X^{1-\delta}}.$$
Using the fact $\theta\in L^1([0,t^*),X^{1-\delta})$ we get
$$\lim_{t<t'\nearrow t^*}\|\theta(t)-\theta(t')\|_{X^{-\delta}}=0,$$
which gives $\theta\in C([0,t^*],X^{-\delta})\cap L^1([0,t^*],X^{1-\delta})$. By applying the first step to the system
$$\begin{cases}
  \partial_t\phi+|D|\phi+u_\phi\nabla \phi  = 0\\
\phi_{/_{t=0}}=\theta(t^*),
   \end{cases}$$
we can extend $\theta$ beyond $t^*$ in $C([0,t^*+\varepsilon],X^{-\delta})\cap L^1([0,t^*+\varepsilon],X^{1-\delta})$, ($\varepsilon>0$ given by step 1),  which contradicts the definition of $t^*$. Then $t^*=+\infty$, and by (\ref{eqt*}), we get
\begin{equation}\label{eqt**}\|\theta(t)\|_{X^{-\delta}}+(1-2^{1-\delta}\|\theta^0\|_{X^0})\int_{0}^{t}\|\theta\|_{X^{1-\delta}}
\leq\|\theta^0\|_{X^{-\delta}},\;\forall t\geq0,\end{equation}
and  $\theta\in C_b(\R^+,X^{-\delta})\cap L^1(\R^+,X^{1-\delta})$.
\subsection{Step3: Asymptotic study in $X^{-\delta}$} In this step we want to prove
\begin{equation}\label{eqL2}\lim_{t\rightarrow\infty}\|\theta(t)\|_{X^{-\delta}}=0.\end{equation}
We do exactly like the first part of the proof of Theorem \ref{th13}: Let $\varepsilon\in(0,1/8)$, then there is a time $t_0$ and a large integer $k$ such that
$$\begin{array}{l}
\displaystyle\int_{0}^\infty\|\theta(t_0+t)\|_{X^1}dt=\int_{t_0}^\infty\|\theta(t)\|_{X^1}dt<\frac{\varepsilon}{4}\\
\displaystyle\int_{A_k^c}|\hat{\theta}(t_0,\xi)|d\xi <\frac{\varepsilon}{2},
\end{array}$$
where
$$A_k=\{\xi\in\mathbb{R}^2/\;|\xi|\leq k\;{\rm and}\;|\hat{\theta}(t_0,\xi)|\leq k\}.$$
We write $\theta(t_0)$ as the sum of two functions as follows $$\theta(t_0)=a^0+b^0$$
where
$$\left\{\begin{array}{lcl}
a^0&=&\mathcal F^{-1}({\bf 1}_{A_k}\widehat{\theta}(t_0))\\
b^0&=&\mathcal F^{-1}({\bf 1}_{A_k^c}\widehat{\theta}(t_0)).\end{array}\right.$$
Clearly, $$\begin{array}{l}
a^0\in L^2\cap X^0\cap X^{-\delta}\\
\|b^0\|_{X^{-\delta}}<\frac{\varepsilon}{2}.
\end{array}$$
Now, consider the following system $$ (S_0)\begin{cases}
  \partial_tb+|D|b+u_b\nabla b  = 0\\
b_{/_{t=0}}=b^0.
   \end{cases}$$
 Combining Theorem \ref{th11} and the above step we get a unique solution $b\in C_b(\mathbb{R}^+,X^0\cap X^{-\delta})\cap L^1(\mathbb{R}^+,X^{1-\delta})$ of  the system $(S_0)$. Moreover, we obtain the following estimation $$\|b\|_{X^{-\delta}}+(1-2^{1-\delta}\|b^0\|_{X^0})\int_{0}^{t}\|b\|_{X^{1-\delta}}\leq\|b^0\|_{X^{-\delta}}.$$
Now put $a(t)=\theta(t_0+t)-b(t)$. As  $\theta,b\in C_b(\mathbb{R}^+,X^0\cap X^{-\delta})\cap L^1(\mathbb{R}^+,X^{1-\delta})$ then  $a\in C_b(\mathbb{R}^+,X^0\cap X^{-\delta})\cap L^1(\mathbb{R}^+,X^{1-\delta})$ and $a$ is the unique solution of the following system
$$\begin{cases}
    \partial a +|D|a +u_a\nabla a+u_a\nabla b+u_b\nabla a& =0 \\
    a_{/_{t=0}} =a^0.
  \end{cases}$$
Using the fact  $a^0\in L^2\cap X^0$ and the second step of the proof of Theorem \ref{th13} (see (\ref{eql22})), we obtain $a\in L^\infty(\R^+,L^2(\R^2))$, and
$$\|a(t)\|_{L^2}\leq 2\|a^0\|_{L^2},\;\forall t\geq0.$$
Define the following subset of $\R^+$, $$F_\varepsilon=\{t\geq 0/\|a(t)\|_{X^{-\delta}}\geq\frac{\varepsilon}{2}\}.$$
Using the interpolation result in Lemma \ref{lm 1.5} with $\sigma=-\delta$, we get
$$\frac{\varepsilon}{2}\leq c\|a(t)\|_{L^2}^{\frac{1}{2-\delta}}\|a(t)\|_{X^{1-\delta}}^{\frac{1-\delta}{2-\delta}},\;\forall t\in F_\varepsilon,$$
and
$$\frac{\varepsilon}{2}\leq 2c\|a^0\|_{L^2}^{\frac{1}{2-\delta}}\|a(t)\|_{X^{1-\delta}}^{\frac{1-\delta}{2-\delta}},\;\forall t\in F_\varepsilon,$$
Then$$(\frac{\varepsilon}{2})^{\frac{2-\delta}{1-\delta}}\leq (2c)^{\frac{2-\delta}{1-\delta}}\|a^0\|_{L^2}^{\frac{1}{1-\delta}}\|a(t)\|_{X^{1-\delta}},\;\forall t\in F_\varepsilon.$$
Integrating over $F_\varepsilon$ we get $$\lambda_1(F_\varepsilon)(\frac{\varepsilon}{2})^{\frac{2-\delta}{1-\delta}}\leq (2c)^{\frac{2-\delta}{1-\delta}}\|a^0\|_{L^2}\|a\|_{L^1(\R^+,X^{1-\delta})}:=M_{\varepsilon,\delta}.$$
Consequently we obtain $$\lambda_1(F_\varepsilon)\leq (\frac{\varepsilon}{2})^{-\frac{2-\delta}{1-\delta}}(2c)^{\frac{2-\delta}{1-\delta}}\|a^0\|_{L^2}\|a\|_{L^1(\R^+,X^{1-\delta})}:=M_{\varepsilon}.$$
As $\lambda_1([0,M_{\varepsilon}+1])=M_{\varepsilon}+1$, there exists $t_\varepsilon\in[0,M_{\varepsilon}+1]\backslash F_\varepsilon$ such that
$\|a(t_\varepsilon)\|_{X^0}< \frac{\varepsilon}{2}.$
Finally \begin{align*}
          \|\theta(t_0+t_\varepsilon)\|_{X^{-\delta}} & \leq \|a(t_\varepsilon)+b(t_\varepsilon)\|_{X^{-\delta}}\\
           & \leq \|a(t_\varepsilon)\|_{X^0}+\|b(t_\varepsilon)\|_{X^{-\delta}}\\
           &  \leq \frac{\varepsilon}{2}+\|b^0\|_{X^{-\delta}}\\
           &  < \frac{\varepsilon}{2}+\frac{\varepsilon}{2}=\varepsilon.
        \end{align*}
Now consider the system starting in $t_0$ $$\begin{cases}
                                            \partial_t\gamma+|D|\gamma+u_\gamma\nabla\gamma=0 \\
                                            \gamma_{/_{t=0}}=\theta(t_0+t_\varepsilon).
                                          \end{cases}$$
By the uniqueness given by Theorem \ref{th11} and the first step (See (\ref{eqt**})), we obtain $\gamma(s)=\theta(t_0+t_\varepsilon+s)$, for $s\geq 0$, and $$\|\gamma(s)\|_{X^{-\delta}}+(1-2^{1-\delta}\|\gamma(0)\|_{X^0})\int_{0}^{t}\|\gamma\|_{X^{1-\delta}}\leq \|\gamma(0)\|_{X^{-\delta}} \;\forall \;s\geq 0.$$
Therefore, the desired result as following $$\|\theta(t)\|_{X^{-\delta}} \leq\|\theta(t_0+t_\varepsilon)\|_{X^{-\delta}} <\varepsilon,\;\forall t\geq t_0+t_\varepsilon.$$

\subsection{Step4: Asymptotic study in $X^0$} Now we want to prove that  $$\|\theta(t)\|_{X^0}= o(t^{-\delta}),\;t\rightarrow\infty.$$
For $\lambda>0$, we have $$\|\theta(t)\|_{X^0} \leq A_\lambda(t)+B_\lambda(t),$$
with $$A_\lambda(t)=\int_{|\xi|<\lambda} e^{-t|\xi|}|\hat{\theta}(t,\xi)|$$
and
$$B_\lambda(t)=\int_{|\xi|<\lambda} e^{-t|\xi|}|\hat{\theta}(t,\xi)|.$$
The low frequencies $A_\lambda(t)$ satisfies
$$A_\lambda\leq \lambda^\delta\|\theta^0\|_{X^{-\delta}}.$$
The analyticity of $\theta$ gives \begin{align*}
      B_\lambda(t)& = \int_{|\xi|>\lambda} e^{-t|\xi|}|\hat{\theta}(t,\xi)|e^{t|\xi|}\\
       & \leq e^{-t\lambda}\int_{|\xi|>\lambda}|\hat{\theta}(t,\xi)|e^{t|\xi|} \\
       & \leq 2\|\theta^0\|_{X^0}e^{-t\lambda}.
    \end{align*}
Consequently we obtain that
\begin{align}
  \|\theta(t)\|_{X^0} & \leq \lambda^\delta\|\theta^0\|_{X^{-\delta}}+ 2\|\theta^0\|_{X^0}e^{-t\lambda}\\
 t^\delta \|\theta(t)\|_{X^0}& \leq (t\lambda)^\delta\|\theta^0\|_{X^{-\delta}}+ 2\|\theta^0\|_{X^0}t^\delta e^{-t\lambda}.
\end{align}
Let $s\geq 0$ and $\phi$ is the solution of the following system
$$\begin{cases}
  \partial_t\phi+|D|\phi+u_\phi\nabla\phi=0 \\
   \phi_{/_{t=0}}=\theta(s).
   \end{cases}$$
Using $(2)$, we get \begin{align*}
                      t^\delta \|\gamma(t)\|_{X^0}& \leq (t\lambda)^\delta\|\gamma(0)\|_{X^{-\delta}}+ 2\|\gamma(0)\|_{X^0}t^\delta e^{-t\lambda} \\
                       t^\delta \|\theta(t+s)\|_{X^0}& \leq (t\lambda)^\delta\|\theta(s)\|_{X^{-\delta}}+ 2\|\theta(s)\|_{X^0}t^\delta e^{-t\lambda}.
                    \end{align*}
Pose that $s=z/2$ and $t=z/2$, we obtain
\begin{align*}
  (z/2)^\delta \|\theta(z)\|_{X^0}& \leq (z/2\lambda)^\delta\|\theta(z/2)\|_{X^{-\delta}}+ 2\|\theta(z/2)\|_{X^0}z^\delta e^{-z\lambda} \\
   \underbrace{z^\delta \|\theta(z)\|_{X^0}}_{=f(z)} & \leq (z\lambda)^\delta\|\theta^0\|_{X^{-\delta}}+ 2\|\theta(z/2)\|_{X^0}(z/2)^\delta e^{-z\lambda}.
\end{align*}
Clearly $f$ is continuous and it satisfies the following inequality
$$f(z)\leq (z\lambda)^\delta\|\theta^0\|_{X^{-\delta}}+ 2e^{-z\lambda}f(z/2).$$
We choose $\lambda=\frac{\log 4}{z},$ to obtain $$f(z)\leq  \underbrace{(\log 4)^\delta\|\theta^0\|_{X^{-\delta}}}_{=A_1}+\frac{1}{2}f(z/2).$$
By Lemma \ref{lm 1.5} we obtain $\sup_{t\geq0}f(t)\leq 2A_1$, and we can deduce that
$$\limsup_{t\rightarrow\infty} t^\delta \|\theta(t)\|_{X^0}\leq2(\log 4)^\delta\|\theta^0\|_{X^{-\delta}}.$$
We do the same work for the system $(S_2)$ with the initial data $\theta(\beta)$, $\beta\geq0$, we get
$$\limsup_{t\rightarrow\infty} t^\delta \|\theta(t+\beta)\|_{X^0}\leq2(\log 4)^\delta\|\theta(\beta)\|_{X^{-\delta}}$$
or
$$\limsup_{t\rightarrow\infty} t^\delta \|\theta(t)\|_{X^0}\leq2(\log 4)^\delta\|\theta(\beta)\|_{X^{-\delta}}.$$
We make $\beta\rightarrow\infty$ and taking into account (\ref{eqL2}), we obtain the requested result.
\section{General remarks}
\subsection{Periodic case} In this section, we give some remarks about the periodic dissipative Quasi-geostrophic equation
$$ (S_3)\begin{cases}
  \partial_t\theta+|D|\theta+u_\theta.\nabla\theta  = 0\;\;{\rm in}\;\mathbb{R}^{+*}\times\mathbb{T}^2\\
u_\theta=\mathcal{R}^\bot\theta=(-R_2\theta,R_1\theta)\\
\theta|_{t=0}=\theta^0,
   \end{cases}$$
   where $\mathbb T^2=(\R/2\pi\Z)^2\simeq[0,2\pi]^2$
In this case we often use $f\in L^1(\mathbb T^2)$ such that
$$\int_{\mathbb T^2}f(x)dx=0.$$
The critical Fourier space $X^0(\mathbb T^2)$ is defined as follows $$X^0(\mathbb T^2)=\{f=\sum_{k\in\Z^2,\,k\neq0}a_k e^{ik.x}\in L^1(\mathbb T^2)/\;\sum_{k\neq0}|a_k|<\infty\}.$$
$X^0(\mathbb T^2)$ is equipped with the norm $$\|f\|_{X^0}=\sum_{k\in\Z^2,\,k\neq0}|\widehat{f}(k)|.$$
We define also$$X^1(\mathbb T^2)=\{f=\sum_{k\in\Z^2,\,k\neq0}a_k e^{ik.x}\in L^1(\mathbb T^2)/\;\sum_{k\neq0}|k|.|a_k|<\infty\}.$$
$X^1(\mathbb T^2)$ is equipped with the norm $$\|f\|_{X^1}=\sum_{k\neq0}|k|.|\widehat{f}(k)|.$$
We use the same methods, we get:\\
{\bf(R1)} If $\|\theta^0\|_{X^0}<1$, then $(S_3)$ has a unique global solution $\theta\in C_b(\R^+,X^0(\mathbb T^2))\cap L^1(\R^+,X^1(\mathbb T^2))$.\\
{\bf(R2)} Moreover, if $\|\theta^0\|_{X^0}<1/4$, then the global solution of $(S_3)$ satisfies
$$\|\theta\|_{\mathfrak{X}^0_{per}}+\|\theta\|_{\mathfrak{X}^1_{per}}\leq 4\|\theta\|_{X^0},$$
where
$$\begin{array}{lcl}
\mathfrak{X}^0_{per}&=&\displaystyle\{((t,x)\mapsto f(t,x))\in \widetilde{L^\infty}(\R^+,X^0(\mathbb T^2))/\,\sum_{k\neq0}\|\widehat{f}(.,k)\|_{L^\infty}<\infty\},\\
\|f\|_{\mathfrak{X}^0_{per}}&=&\displaystyle\sum_{k\neq0}\|\widehat{f}(.,k)\|_{L^\infty},\\
\mathfrak{X}^1_{per}&=&\displaystyle L^1(\R^+,X^1(\mathbb T^2)),\\
\|f\|_{\mathfrak{X}^1_{per}}&=&\displaystyle\int_{0}^\infty\sum_{k\neq0}|\widehat{f}(t,k)|dt.
\end{array}$$
{\bf(R3)}  If $\|\theta^0\|_{X^0}<1$, then the global solution of $(S_3)$ satisfies
$$\lim_{t\rightarrow\infty}\|\theta(t)\|_{X^0}=0.$$
{\bf(R4)}  If $\|\theta^0\|_{X^0}<1/8$, then the global solution of $(S_3)$ satisfies
$$\|e^{\frac{t}{2}|D|}\theta\|_{\mathfrak{X}^0_{per}}+\|e^{\frac{t}{2}|D|}\theta\|_{\mathfrak{X}^1_{per}}\leq 12\|\theta\|_{X^0}.$$
Particulary, we get $e^{\frac{t}{2}}\|\theta(t)\|_{X^0}\leq 12\|\theta\|_{X^0}$ and $\|\theta(t)\|_{X^0}\leq 12\|\theta\|_{X^0}e^{-\frac{t}{2}}$. This property of exponential decay is due to the fact that $\mathbb T^2$ is compact and zero is isolated in the frequencies space. We can improve the analytical index $1/2$ to $r\in(0,1)$: There is $\varepsilon_r\in(0,1)$ such that if $\|\theta^0\|_{X^0}<\varepsilon_r$, then the global solution $\theta$ of $(S_3)$ satisfies
$$\|e^{\varepsilon_rt|D|}\theta\|_{\mathfrak{X}^0_{per}}+\|e^{\varepsilon_rt|D|}\theta\|_{\mathfrak{X}^1_{per}}\leq C_r\|\theta\|_{X^0},$$
which implies $\|\theta(t)\|_{X^0}\leq C_r\|\theta^0\|_{X^0}e^{-rt}$.\\
{\bf(R5)} Combining {\bf(R3)} and {\bf(R4)}, we obtain: If $\|\theta^0\|_{X^0}<1/4$, then there is a time $t_0\geq 0$ such that $\|\theta(t_0)\|_{X^0}<1/8$. Moreover the global solution of $(S_3)$ satisfies
$$\|e^{\frac{t}{2}|D|}\theta(t_0+.)\|_{\mathfrak{X}^0_{per}}+\|e^{\frac{t}{2}|D|}\theta(t_0+.)\|_{\mathfrak{X}^1_{per}}\leq 12\|\theta(t_0)\|_{X^0}.$$
Similarly, if $r\in(0,1)$ there is a time $t_r\geq0$ such that $\|\theta(t_r)\|_{X^0}<\varepsilon_r$ and $\theta$ satisfies
$$\|e^{\varepsilon_rt|D|}\theta(t_r+.)\|_{\mathfrak{X}^0_{per}}+\|e^{\varepsilon_rt|D|}\theta(t_r+.)\|_{\mathfrak{X}^1_{per}}\leq C_r\|\theta(t_r)\|_{X^0},$$
which implies $\|\theta(t_r+t)\|_{X^0}\leq C_r\|\theta(t_r)\|_{X^0}e^{-rt}$.
\subsection{Property of analyticity } Generally, the property of analyticity comes from the linear part of the nonlinear equation
$$ (S)\begin{cases}
  \partial_tf+|D|^{\sigma}f+u_f.\nabla f= 0\;\;{\rm in}\;\mathbb{R}^{+*}\times\Omega \\
f|_{t=0}=f_0,
   \end{cases}$$
where $\Omega=\R^2$ or $\mathbb T^2$. If $f_0\in X^0(\Omega)$ (rep. $f_0\in H^1(\Omega)$),  then $g(t)=e^{-t|D|^{\sigma}}f_0$,\;$\sigma\geq1$, is the unique solution of the following system
$$ (LS)\begin{cases}
  \partial_tg+|D|^{\sigma}g= 0\;\;{\rm in}\;\mathbb{R}^{+*}\times\Omega \\
g|_{t=0}=f_0.
   \end{cases}$$
Then $e^{t|D|^{\sigma}}g(t)=f_0\in C_b(\R^+,X^0)$ (rep. $C_b(\R^+,H^1)$). But for the system $(S)$ even if $\sigma>1$ we don't hope to have more than the analyticity of his solution. This is a technical problem of analysis: For $p>0$, the property
$$|\xi|^p\leq |\xi-\eta|^p+|\eta|^p,\;\forall \xi,\eta\in\R^2$$
or
$$|k|^p\leq |k-n|^p+|n|^p,\;\forall k,n\in\mathbb Z^2$$
implies $p\in(0,1]$.\\
Indeed: Suppose that $p>1$. For $\beta,\gamma\in\N$ and by applying the above inequality with $k=(\beta+\gamma,0)$ and $n=(\gamma,0)$ we get
$$(\beta+\gamma)^p\leq \beta^p+\gamma^p$$
and
$$(1+\frac{\gamma}{\beta})^p\leq 1+(\frac{\gamma}{\beta})^p.$$
By density of $\mathbb Q$ in $\R$, we obtain
$$(1+a)^p\leq 1+a^p,\;\forall a\in(0,\infty),$$
which implies
$$\frac{(1+a)^p-1}{a}\leq a^{p-1},\;\forall a\in(0,\infty).$$
By taking the limit $a\rightarrow0^+$, we obtain contradiction.\\
In conclusion with classical analysis, one cannot prove the following type of inequality, if $p>1$ and $z>0$,
$$\|e^{z|D|^p}(fg)\|_{X^0}\leq C_p\|e^{z|D|^p}f\|_{X^0}\|e^{z|D|^p}g\|_{X^0}.$$
At this stage, we recall the result of Foias and Temam (See \cite{FT1}) with $\sigma = 2$ for the Navier-Stokes equations on $\mathbb T^3$, where they proved only the analyticity of the solution. So, in our case we don't hope to have more than analyticity of the global solution.

\end{document}